# PARTIALLY OBSERVED INFORMATION AND INFERENCE ABOUT NON-GAUSSIAN MIXED LINEAR MODELS[1]

### By Jiming Jiang

*University of California, Davis*

In mixed linear models with nonnormal data, the Gaussian Fisher information matrix is called a quasi-information matrix (QUIM). The QUIM plays an important role in evaluating the asymptotic covariance matrix of the estimators of the model parameters, including the variance components. Traditionally, there are two ways to estimate the information matrix: the estimated information matrix and the observed one. Because the analytic form of the QUIM involves parameters other than the variance components, for example, the third and fourth moments of the random effects, the estimated QUIM is not available. On the other hand, because of the dependence and nonnormality of the data, the observed QUIM is inconsistent. We propose an estimator of the QUIM that consists partially of an observed form and partially of an estimated one. We show that this estimator is consistent and computationally very easy to operate. The method is used to derive large sample tests of statistical hypotheses that involve the variance components in a non-Gaussian mixed linear model. Finite sample performance of the test is studied by simulations and compared with the delete-group jackknife method that applies to a special case of non-Gaussian mixed linear models.

**1. Introduction.** Mixed linear models are widely used in practice, especially in situations involving correlated observations. A typical assumption regarding these models is that the observations are normally distributed, or, equivalently, that the random effects and errors in the model are normal. However, as is well known, the normality assumption is likely to be violated. For example, Lange and Ryan [19] gave several examples that show that nonnormality of the random effects is, indeed, encountered in practice. The authors also developed a method for assessing normality of the random effects.

Received July 2003; revised January 2005.

[1]Supported in part by NSF Grant DMS-02-03676.

Due to such concerns, some researchers have considered the use of Gaussian maximum likelihood (ML) or restricted maximum likelihood (REML) estimators in nonnormal situations; see Richardson and Welsh [22], Jiang [12, 13] and Heyde [10, 11], among others. Throughout this paper these estimators will be called ML and REML estimators even if normality does not hold. In particular, Jiang [12, 13] established consistency and asymptotic normality of REML and ML estimators in nonnormal situations under regularity conditions. Furthermore, Jiang [12] derived the asymptotic covariance matrix (ACM) of the REML estimator of the variance components as well as that of the ML estimator without assuming normality. Also see [14]. The ACM is important for various inferences about the model parameters, including interval estimation and hypothesis testing. Unfortunately, the ACM under nonnormality involves parameters other than the variance components, for example, the third and fourth moments of the random effects. Note that standard procedures such as ML and REML do not produce estimators of these additional parameters. For years this complication has undermined the potential usefulness of the ACM in nonnormal situations.

To see exactly where the problem occurs, consider the mixed linear model

$$(1) \qquad y = X\beta + Z_1\alpha_1 + \cdots + Z_s\alpha_s + \varepsilon,$$

where $y$ is an $N \times 1$ vector of observations, $X, Z_1, \ldots, Z_s$ are known matrices, $\beta$ is a $p \times 1$ vector of unknown parameters (the fixed effects), $\alpha_1, \ldots, \alpha_s$ are vectors of random effects and $\varepsilon$ is a vector of errors. It is assumed that $\alpha_1, \ldots, \alpha_s, \varepsilon$ are independent. Furthermore, the components of $\alpha_j$ are i.i.d. with mean 0 and variance $\sigma_j^2$, $1 \le j \le s$, and the components of $\varepsilon$ are i.i.d. with mean 0 and variance $\sigma_0^2$. Without loss of generality, let $\mathrm{rank}(X) = p$. Note that normality is not assumed in this model, nor is any other specific distribution assumed. Also, w.l.o.g. consider the Hartley–Rao form of the variance components [9]: $\lambda = \sigma_0^2$, $\gamma_j = \sigma_j^2/\sigma_0^2$, $1 \le j \le s$. Let $\theta = (\lambda, \gamma_1, \ldots, \gamma_s)'$.

According to [12], the ACM of the REML estimator $\hat{\theta}$ is given by

$$(2) \qquad \Sigma_{\mathrm{R}} = \left\{ \mathrm{E}\left( \frac{\partial^2 l_{\mathrm{R}}}{\partial\theta\,\partial\theta'} \right) \right\}^{-1} \mathrm{Var}\left( \frac{\partial l_{\mathrm{R}}}{\partial\theta} \right) \left\{ \mathrm{E}\left( \frac{\partial^2 l_{\mathrm{R}}}{\partial\theta\,\partial\theta'} \right) \right\}^{-1},$$

where $l_{\mathrm{R}}$ is the Gaussian restricted log-likelihood function, that is, the log-likelihood based on $z = T'y$, where $y$ satisfies (1) with normally distributed random effects and errors, and $T$ is an $N \times (N - p)$ matrix of full rank such that $T'X = 0$. The matrix $\mathcal{I}_2 = \mathrm{E}(\partial^2 l_{\mathrm{R}}/\partial\theta\,\partial\theta')$ depends only on $\theta$, whose estimator is already available. However, unlike $\mathcal{I}_2$, the matrix $\mathcal{I}_1 = \mathrm{Var}(\partial l_{\mathrm{R}}/\partial\theta)$ depends on, in addition to $\theta$, the kurtoses of the random effects and errors. Similarly, let $\psi = (\beta'\theta')'$ and let $\hat{\psi}$ be the ML estimator of $\psi$. By the result of Jiang [14], it can be shown that the ACM of $\hat{\psi}$ is given by

$$\Sigma = \left\{ E\left( \frac{\partial^2 l}{\partial \psi \, \partial \psi'} \right) \right\}^{-1} \text{Var}\left( \frac{\partial l}{\partial \psi} \right) \left\{ E\left( \frac{\partial^2 l}{\partial \psi \, \partial \psi'} \right) \right\}^{-1},$$

where $l$ is the Gaussian log-likelihood. Here, again, the matrix $\mathcal{I}_2 = E(\partial^2 l / \partial \psi \, \partial \psi')$ depends only on $\theta$, but the matrix $\mathcal{I}_1 = \text{Var}(\partial l / \partial \psi)$ depends on, in addition to $\theta$, the kurtoses as well as the third moments of random effects and errors.

It is clear that the key issue is how to estimate $\mathcal{I}_1$, which we call the quasi-information matrix (QUIM) for an obvious reason. Consider, for example, the ML case. If $l$ were the true log-likelihood, then we would have $\mathcal{I}_1 = -\mathcal{I}_2$, which is the Fisher information matrix. Traditionally, there are two ways to estimate the Fisher information: (i) the estimated information and (ii) the observed information. See, for example, [7] for a discussion and comparison of the two methods in the i.i.d. case. It is known that standard procedures in mixed model analysis such as ML and REML do not produce estimators of the third and fourth moments of the random effects and errors. Therefore, according to our previous discussion, method (i) is not possible unless one finds some way to estimate these higher moments. Assuming that the random effects and errors are symmetrically distributed, in which case the third moments vanish, Jiang [15] proposed an empirical method of moments (EMM) to estimate the kurtoses of the random effects and errors. It is clear that this method has a limitation, because, like normality, symmetry may not hold in practice. When the third moments are nonzero, the EMM cannot be used. Furthermore, the situation to which the EMM applies is somewhat restrictive and requires certain orthogonal decompositions of the linear spaces generated by the design matrices of the random effects. Simulation results have suggested that the EMM estimator may have large variance even when the sample size is moderately large. As for method (ii), it is not all that clear how this should be defined in cases of correlated observations. For simplicity, let us assume that $\psi$ is a scalar. With independent observations, we have

$$\mathcal{I}_1 = E\left\{ \sum_{i=1}^{N} \left( \frac{\partial l_i}{\partial \psi} \right)^2 \right\}, \tag{3}$$

where $l_i$ is the log-likelihood based on $y_i$, the $i$th observation. Therefore, an observed information is $\tilde{\mathcal{I}}_1 = \sum_{i=1}^{N} (\partial l_i / \partial \psi |_{\hat{\psi}})^2$, where $\hat{\psi}$ is the ML estimator. This is a consistent estimator of $\mathcal{I}_1$ in the sense that $\tilde{\mathcal{I}}_1 - \mathcal{I}_1 = o_P(\mathcal{I}_1)$ or, equivalently, $\tilde{\mathcal{I}}_1 / \mathcal{I}_1 \to 1$ in probability. However, if the observations are correlated, (3) does not hold. In this case, since $\mathcal{I}_1 = E\{(\partial l / \partial \psi)^2\}$, one might attempt to define $\tilde{\mathcal{I}}_1 = (\partial l / \partial \psi |_{\hat{\psi}})^2$. However, this is zero since $\hat{\psi}$ is the MLE. Even if $\hat{\psi}$ is a different (consistent) estimator, the expression is not a consistent estimator. For example, in the independent case this is the same as $(\sum_{i=1}^{N} \partial l_i / \partial \psi|_{\hat{\psi}})^2$, which, asymptotically, is equivalent to $N$ times the square of a normal random variable. Therefore, it is not true that $\hat{\mathcal{I}}_1 - \mathcal{I}_1 = o_P(\mathcal{I}_1)$. Alternatively, if normality holds, one may define $l_i$ as the logarithm of the conditional density of $y_i$ given $y_1, \ldots, y_{i-1}$. It follows that $\partial l_i / \partial \psi$, $1 \le i \le N$, is a sequence of martingale differences [with respect to the $\sigma$-fields $\mathcal{F}_i = \sigma(y_1, \ldots, y_i)$, $1 \le i \le N$]. Thus, we still have (3) but with new definitions of $l_i$'s; hence $\hat{\mathcal{I}}_1$ can be defined similarly as in the independent case. However, if normality does not hold, this latter strategy also does not work (because $\partial l_i / \partial \psi$ is no longer a martingale difference).

We now explain our approach to the problem using a simple example.

EXAMPLE 1. Consider the following model with crossed random effects: $y_{ij} = \mu + v_i + w_j + e_{ij}$, $i = 1, \ldots, m$, $j = 1, \ldots, n$, where $\mu$ is an unknown mean, $v_i$ and $w_j$ are random effects, and $e_{ij}$ is an error. It is assumed that the $v_i$'s are i.i.d. with mean 0 and variance $\sigma_1^2$, the $w_j$'s are i.i.d. with mean 0 and variance $\sigma_2^2$, the $e_{ij}$'s are i.i.d. with mean 0 and variance $\sigma_0^2$, and $v$, $w$ and $e$ are independent. Consider an element of the QUIM $\text{Var}(\partial l_R / \partial \theta)$ for REML estimation, say, $\text{var}(\partial l_R / \partial \lambda)$, where $l_R$ is the Gaussian restricted log-likelihood and $\theta = (\lambda, \gamma_1, \gamma_2)'$ ($\lambda$ and $\gamma$'s as defined earlier). By the result of Jiang ([16], Section 5, Example 2), it can be shown that $\partial l_R / \partial \lambda = \{u'Bu - (mn-1)\lambda\}/2\lambda^2$, where $u = y - \mu \mathbf{1}_m \otimes \mathbf{1}_n$ with $y = (y_{ij})_{1 \le i \le m, 1 \le j \le n}$ (as a vector in which the components are ordered as $y_{11}, \ldots, y_{1n}, y_{21}, \ldots$) and

$$B = I_m \otimes I_n - \frac{1}{n}\left(1 - \frac{1}{1 + \gamma_1 n}\right) I_m \otimes J_n - \frac{1}{m}\left(1 - \frac{1}{1 + \gamma_2 m}\right) J_m \otimes I_n$$
$$+ \frac{1}{mn}\left(1 - \frac{1}{1 + \gamma_1 n} - \frac{1}{1 + \gamma_2 m}\right) J_m \otimes J_n$$
$$= I_m \otimes I_n + \lambda_1 I_m \otimes J_n + \lambda_2 J_m \otimes I_n + \lambda_3 J_m \otimes J_n.$$

Hereafter, $I_n$ and $\mathbf{1}_n$ represent the $n$-dimensional identity matrix and the vector of 1's, respectively, $J_n = \mathbf{1}_n \mathbf{1}_n'$ and $\otimes$ means Kronecker product. Define $\kappa_0 = \text{E}(e_{11}^4) - 3\lambda^2$, $\kappa_1 = \text{E}(v_1^4) - 3\lambda^2 \gamma_1^2$, $\kappa_2 = \text{E}(w_1^4) - 3\lambda^2 \gamma_2^2$ (note that these are the kurtoses) and $t_0 = 1 + \lambda_1 + \lambda_2 + \lambda_3$, $t_1 = \{(m-1)n\}/\{m(1 + \gamma_1 n)\}$, $t_2 = \{m(n-1)\}/n(1 + \gamma_2 m)\}$; $m_0 = mn$, $m_1 = m$ and $m_2 = n$. By Lemma 1 in the sequel, it can be shown that

$$\text{var}\left(\frac{\partial l_R}{\partial \lambda}\right)$$
$$= \text{E}\left\{(a_0 + a_1 + a_2)\sum_{i,j} u_{ij}^4 - a_1 \sum_i \left(\sum_j u_{ij}\right)^4 - a_2 \sum_j \left(\sum_i u_{ij}\right)^4\right\}$$

(4)
$$+ \left[ \frac{mn-1}{2\lambda^2} - \frac{3mnt_0^2}{4\lambda^2} \{ (1 + \gamma_1 + \gamma_2)^2 - (t_3 + t_4) \} - \frac{3(t_1^2 t_3 m + t_2^2 t_4 n)}{4\lambda^2} \right]$$

$$= S_1 + S_2,$$

where

$$a_0 = \frac{t_0^2}{4\lambda^4}, \qquad a_1 = \frac{nt_0^2 - t_1^2}{4\lambda^4 n(n^3-1)}, \qquad a_2 = \frac{mt_0^2 - t_2^2}{4\lambda^4 m(m^3-1)},$$

$$t_3 = \frac{n(1 + \gamma_2 + \gamma_1 n)^2 - (1 + \gamma_1 + \gamma_2)^2}{n^3 - 1},$$

$$t_4 = \frac{m(1 + \gamma_1 + \gamma_2 m)^2 - (1 + \gamma_1 + \gamma_2)^2}{m^3 - 1}.$$

It is clear that $S_2$ can be estimated by replacing the variance components by their REML estimators, which are already available. As for $S_1$, it cannot be estimated in the same way for the reason given above. However, the form of $S_1$ [cf. with (3)] suggests an "observed" estimator by taking out the expectation sign and replacing the parameters involved by their REML estimators. In fact, as $m, n \to \infty$, this observed $S_1$, say, $\hat{S}_1$, is consistent in the sense that $\hat{S}_1/S_1 \to 1$ in probability. It is interesting to note that $S_2$ cannot be consistently estimated by an observed form. In conclusion, $S_1$ cannot be estimated by an estimated form, but can be estimated by an observed form; $S_2$ can be estimated by an estimated form, but not by an observed form. Thus, we have reached a balance.

We propose to use such a method to estimate the QUIM. Because the estimator consists partially of an observed form and partially of an estimated one, it is called a partially observed quasi-information matrix (POQUIM).

One application of POQUIM is to derive robust dispersion tests in mixed linear models. A dispersion test is a test of a statistical hypothesis that involves the variance components. Such tests, exact or asymptotic, are available in the literature (e.g., [17, 23]), but only under the normality assumption. Since the latter is likely to be violated in practice, as a robust approach, it is of interest to derive dispersion tests that do not rely on normality. Using the results of Jiang [12, 14], it is possible to derive an asymptotic dispersion test based on either the REML or the ML estimators without assuming normality, provided that the ACM can be consistently estimated. The POQUIM will provide such a consistent estimator.

The rest of the paper is organized as follows. In Section 2 we explain how one comes up with the decomposition (4), that is, we derive POQUIM for a general non-Gaussian mixed linear model with REML estimation of the variance components. Sufficient conditions will be given for the consistency of POQUIM as well as an estimator of the ACM of the REML estimator. In Section 3 we use several examples to illustrate the main results of Section 2. In Section 4 we consider POQUIM for ML estimation. In Section 5 we apply POQUIM to robust dispersion tests in mixed linear models. Some simulated examples are considered in Section 6, in which we study the finite sample performance of POQUIM in the context of robust dispersion tests and compare it with the delete-group jackknife method of Arvesen [1] (also see [2]) in a case where the latter applies. In Section 7 we discuss extension of POQUIM to quasi-likelihood estimation and remark on other issues. Proofs and other technical details are given in Section 8.

## 2. POQUIM for REML.

The REML case is relatively simple compared to ML, because only estimation of the variance components is involved. Furthermore, as will be seen, the QUIM in this case does not involve the third moments of the random effects and errors.

Under model (1) and normality, the restricted log-likelihood for estimating the variance components $\lambda$ and $\gamma_j$, $1 \le j \le s$, is

$$l_{\mathrm{R}}(\theta) = c - \tfrac{1}{2}\{\log(|T'VT|) + y'Py\},$$

where $\theta = (\lambda, \gamma_1, \ldots, \gamma_s)'$, $c$ is a constant, $V = \mathrm{Var}(y) = \lambda(I + \gamma_1 Z_1 Z_1' + \cdots + \gamma_s Z_s Z_s')$ ($I$ is the $N \times N$ identity matrix), $T$ is any $N \times (N-p)$ matrix such that $\mathrm{rank}(T) = N-p$ and $T'X = 0$ ($|\cdot|$ means determinant), and $P = T(T'VT)^{-1}T' = V^{-1} - V^{-1}X(X'V^{-1}X)^{-1}X'V^{-1}$ (e.g., [23], page 451). If normality does not hold, (5) is not the true restricted log-likelihood but, instead, the quasi-restricted log-likelihood. It is shown in Section 8.1 that $\partial l_{\mathrm{R}}/\partial \theta_j = u'B_j u - b_j$, $0 \le j \le s$, where $\theta_0 = \lambda$, $\theta_j = \gamma_j$, $1 \le j \le s$; $u = y - X\beta$; $B_0 = (2\lambda)^{-1}P$, $B_j = (\lambda/2)PZ_jZ_j'P$, $1 \le j \le s$; $b_0 = (N-p)/2\lambda$ and $b_j = (\lambda/2)\mathrm{tr}(PZ_jZ_j')$, $1 \le j \le s$. Note that $b_j = \mathrm{E}(u'B_j u)$, $0 \le j \le s$.

### 2.1. Derivation.

Let $u_i = y_i - x_i'\beta$ be the $i$th component of $u$, where $x_i'$ is the $i$th row of $X$. The kurtoses of the random effects and errors are defined as $\kappa_t = \mathrm{E}(\alpha_{t1}^4) - 3\sigma_t^4 = \mathrm{E}(\alpha_{t1}^4) - 3(\lambda\gamma_t)^2$, $0 \le t \le s$, where $\alpha_0 = \varepsilon$ and $\gamma_0 = 1$. Also, with a slight abuse of the notation, let $z_{it}'$ and $z_{tl}$ be the $i$th row and $l$th column of $Z_t$, respectively, $0 \le t \le s$, where $Z_0 = I$. Define $\Gamma(i_1, i_2) = \sum_{t=0}^{s} \gamma_t(z_{i_1 t} \cdot z_{i_2 t})$. Here the dot product of vectors $a_1, \ldots, a_k$ of the same dimension is defined as $a_1 \cdot a_2 \cdots a_k = \sum_l a_{1l} a_{2l} \cdots a_{kl}$. Also, let $m_t$ be the dimension of $\alpha_t$, $0 \le t \le s$ (so that $m_0 = N$). We begin with an expression for $\mathrm{cov}(u_{i_1} u_{i_2}, u_{i_3} u_{i_4})$ ($1 \le i_1, \ldots, i_4 \le N$) as well as one for $\mathrm{cov}(\partial l_{\mathrm{R}}/\partial \theta_j, \partial l_{\mathrm{R}}/\partial \theta_k)$, the $(j, k)$ element of $\mathcal{I}_1$.

LEMMA 1. *We have*

$$\mathrm{cov}(u_{i_1} u_{i_2}, u_{i_3} u_{i_4})$$
$$= \lambda^2\{\Gamma(i_1, i_3)\Gamma(i_2, i_4) + \Gamma(i_1, i_4)\Gamma(i_2, i_3)\} + \sum_{t=0}^{s}\kappa_t z_{i_1 t} \cdots z_{i_4 t},$$

where $z_{i_1 t} \cdots z_{i_4 t} = z_{i_1 t} \cdot z_{i_2 t} \cdot z_{i_3 t} \cdot z_{i_4 t}$. Furthermore, we have

$$(7) \quad \mathrm{cov}\left(\frac{\partial l_{\mathrm{R}}}{\partial \theta_j}, \frac{\partial l_{\mathrm{R}}}{\partial \theta_k}\right) = 2\,\mathrm{tr}(B_j V B_k V) + \sum_{t=0}^{s} \kappa_t \sum_{l=1}^{m_t} (z'_{tl} B_j z_{tl})(z'_{tl} B_k z_{tl}).$$

The proof is given in Section 8.2.

Let $f_1, \ldots, f_L$ be the different nonzero functional values of

$$(8) \quad f(i_1, \ldots, i_4) = \sum_{t=0}^{s} \kappa_t z_{i_1 t} \cdots z_{i_4 t}.$$

Note that this is the second term on the right-hand side of (6). Here functional value means $f(i_1, \ldots, i_4)$ as a function of $\kappa = (\kappa_t)_{0 \le t \le s}$. For example, $\kappa_0 + \kappa_1$ and $\kappa_2 + \kappa_3$ are different functions (even if their values may be the same for some $\kappa$). Also, let 0 denote the zero function (of $\kappa$). Then without using (7) we have

$$
\begin{aligned}
\mathrm{cov}\left(\frac{\partial l_{\mathrm{R}}}{\partial \theta_j}, \frac{\partial l_{\mathrm{R}}}{\partial \theta_k}\right) &= \sum_{i_1, \ldots, i_4} B_{j, i_1, i_2} B_{k, i_3, i_4}\,\mathrm{cov}(u_{i_1} u_{i_2}, u_{i_3} u_{i_4}) \\
(9) \qquad &= \sum_{f(i_1, \ldots, i_4) = 0} B_{j, i_1, i_2} B_{k, i_3, i_4}\,\mathrm{cov}(u_{i_1} u_{i_2}, u_{i_3} u_{i_4}) \\
&\quad + \sum_{l=1}^{L} \sum_{f(i_1, \ldots, i_4) = f_l} B_{j, i_1, i_2} B_{k, i_3, i_4}\,\mathrm{cov}(u_{i_1} u_{i_2}, u_{i_3} u_{i_4}) \\
&= \sum_{l=0}^{L} S_l
\end{aligned}
$$

with $S_l$, $0 \le l \le L$, defined in obvious ways. According to Lemma 1, the left-hand side of (9) depends on the higher moments only through $\kappa$. By (6) and (8) we have

$$(10) \quad S_0 = 2\lambda^2 \sum_{f(i_1, \ldots, i_4) = 0} B_{j, i_1, i_2} B_{k, i_3, i_4} \Gamma(i_1, i_3) \Gamma(i_2, i_4),$$

which depends only on $\theta$. Furthermore, for $1 \le l \le L$ write

$$
\begin{aligned}
S_l &= c_l \sum_{f(i_1, \ldots, i_4) = f_l} \mathrm{cov}(u_{i_1} u_{i_2}, u_{i_3} u_{i_4}) \\
&\quad + \sum_{f(i_1, \ldots, i_4) = f_l} (B_{j, i_1, i_2} B_{k, i_3, i_4} - c_l)\,\mathrm{cov}(u_{i_1} u_{i_2}, u_{i_3} u_{i_4}) \\
&= S_{l,1} + S_{l,2},
\end{aligned}
$$

where $c_l$ is a constant to be determined later on. By (6) we have

$$S_{l,2} = \sum_{f(i_1,\ldots,i_4)=f_l} (B_{j,i_1,i_2}B_{k,i_3,i_4} - c_l)[f_l + \lambda^2\{\cdots\}]$$

$$= f_l \sum_{f(i_1,\ldots,i_4)=f_l} (B_{j,i_1,i_2}B_{k,i_3,i_4} - c_l) + \cdots,$$

where $\cdots$ depends only on $\theta$. If we let the coefficient of $f_l$ in the above be equal to zero, we have

$$(11) \qquad c_l = \frac{1}{|\{f(i_1,\ldots,i_4)=f_l\}|} \sum_{f(i_1,\ldots,i_4)=f_l} B_{j,i_1,i_2}B_{k,i_3,i_4},$$

where $|\cdot|$ denotes cardinality. With this choice of $c_l$, we have

$$S_{l,2} = \lambda^2 \sum_{f(i_1,\ldots,i_4)=f_l} (B_{j,i_1,i_2}B_{k,i_3,i_4} - c_l)\{\Gamma(i_1,i_3)\Gamma(i_2,i_4) + \Gamma(i_1,i_4)\Gamma(i_2,i_3)\}$$

$$= 2\lambda^2 \sum_{f(i_1,\ldots,i_4)=f_l} (B_{j,i_1,i_2}B_{k,i_3,i_4} - c_l)\Gamma(i_1,i_3)\Gamma(i_2,i_4),$$

which depends only on $\theta$. Note that $c_l$ depends only on $\theta$. On the other hand, by the fact that $\mathrm{E}(u_{i_1}u_{i_2}) = \lambda\Gamma(i_1,i_2)$ (see the first paragraph of Section 8.2), we have

$$S_{l,1} = c_l \sum_{f(i_1,\ldots,i_4)=f_l} \{\mathrm{E}(u_{i_1}\cdots u_{i_4}) - \lambda^2\Gamma(i_1,i_2)\Gamma(i_3,i_4)\}$$

$$= \mathrm{E}\left(c_l \sum_{f(i_1,\ldots,i_4)=f_l} u_{i_1}\cdots u_{i_4}\right) - \lambda^2 c_l \sum_{f(i_1,\ldots,i_4)=f_l} \Gamma(i_1,i_3)\Gamma(i_2,i_4).$$

Note that $\sum_{f(i_1,\ldots,i_4)=f_l} \Gamma(i_1,i_2)\Gamma(i_3,i_4) = \sum_{f(i_1,\ldots,i_4)=f_l} \Gamma(i_1,i_3)\Gamma(i_2,i_4)$, because $f(i_1,\ldots,i_4)$ is symmetric in $i_1,\ldots,i_4$. Therefore, we have, by combining the above,

$$\begin{aligned}(12) \qquad S_l = {}&\mathrm{E}\left(c_l \sum_{f(i_1,\ldots,i_4)=f_l} u_{i_1}\cdots u_{i_4}\right) \\ &+ 2\lambda^2 \sum_{f(i_1,\ldots,i_4)=f_l} B_{j,i_1,i_2}B_{k,i_3,i_4}\Gamma(i_1,i_3)\Gamma(i_2,i_4) \\ &- 3\lambda^2 c_l \sum_{f(i_1,\ldots,i_4)=f_l} \Gamma(i_1,i_3)\Gamma(i_2,i_4).\end{aligned}$$

Note that $c_l$ defined by (11) depends on $j$ and $k$, that is, $c_l = c_{j,k,l}$. If we define $c_{j,k}(i_1,\ldots,i_4) = c_{j,k,l}$, if $f(i_1,\ldots,i_4) = f_l$, $1 \le l \le L$, then by (9),

(10) and (12) it can be shown that

$$\text{cov}\left(\frac{\partial l_{\mathrm{R}}}{\partial \theta_j}, \frac{\partial l_{\mathrm{R}}}{\partial \theta_k}\right) = \mathrm{E}\left\{\sum_{f(i_1,\ldots,i_4)\neq 0} c_{j,k}(i_1,\ldots,i_4)u_{i_1}\cdots u_{i_4}\right\} + 2\,\text{tr}(B_j V B_k V)$$
$$- 3\lambda^2 \sum_{f(i_1,\ldots,i_4)\neq 0} c_{j,k}(i_1,\ldots,i_4)\Gamma(i_1,i_3)\Gamma(i_2,i_4).$$

We summarize the result in terms of a theorem. Write $\mathcal{I}_{1,jk} = \text{cov}(\partial l_{\mathrm{R}}/\partial \theta_j, \partial l_{\mathrm{R}}/\partial \theta_k)$, which is the $j,k$ element of the QUIM $\mathcal{I}_1 = \text{Var}(\partial l_{\mathrm{R}}/\partial \theta)$.

THEOREM 1. *For any non-Gaussian mixed linear model* (1), *we have*

$$
\begin{aligned}
\mathcal{I}_{1,jk} &= 2\,\text{tr}(B_j V B_k V) + \sum_{t=0}^{s}\kappa_t\sum_{l=1}^{m_t}(z'_{tl}B_j z_{tl})(z'_{tl}B_k z_{tl}) \\
&= \mathrm{E}\left\{\sum_{f(i_1,\ldots,i_4)\neq 0} c_{j,k}(i_1,\ldots,i_4)u_{i_1}\cdots u_{i_4}\right\} \\
&\quad + \left\{2\,\text{tr}(B_j V B_k V) - 3\lambda^2 \sum_{f(i_1,\ldots,i_4)\neq 0} c_{j,k}(i_1,\ldots,i_4)\Gamma(i_1,i_3)\Gamma(i_2,i_4)\right\} \\
&= \mathcal{I}_{1,1,jk} + \mathcal{I}_{1,2,jk},
\end{aligned}
$$

(13)

$0 \leq j,k \leq s$, *where* $c_{j,k}(i_1,\ldots,i_4) = c_{j,k,l}$, *if* $f(i_1,\ldots,i_4) = f_l$, $1 \leq l \leq L$, *with*

$$(14) \qquad c_{j,k,l} = \frac{1}{|\{f(i_1,\ldots,i_4) = f_l\}|}\sum_{f(i_1,\ldots,i_4)=f_l} B_{j,i_1,i_2}B_{k,i_3,i_4}.$$

Of course, (13) can be verified directly, but the derivation above also explains where the thought came from. Note that $2\,\text{tr}(B_j V B_k V)$ is the Gaussian covariance between $\partial l_{\mathrm{R}}/\partial \theta_j$ and $\partial l_{\mathrm{R}}/\partial \theta_k$. This means that under normality $\mathcal{I}_{1,1,jk}$ is identical to the second term in $\mathcal{I}_{1,2,jk}$ with the negative sign removed. Of course, this can be easily verified using (6). On the other hand, without normality $\mathcal{I}_{1,1,jk}$ may involve higher moments of the random effects and errors, and this is why the expectation is not taken inside the summation. Instead, we propose to estimate $\mathcal{I}_{1,1,jk}$ by taking out the expectation sign and replacing any parameter involved by its REML estimator, that is, $\hat{\mathcal{I}}_{1,1,jk} = \sum_{f(i_1,\ldots,i_4)\neq 0}\hat{c}_{j,k}(i_1,\ldots,i_4)\hat{u}_{i_1}\cdots\hat{u}_{i_4}$, where $\hat{c}_{j,k}(i_1,\ldots,i_4)$ is defined in the same way as $c_{j,k}(i_1,\ldots,i_4)$ except with $\theta$ replaced by $\hat{\theta}$, and $\hat{u}_i = y_i - x'_i\hat{\beta}$. Here $\hat{\theta}$ is the REML estimator of $\theta$, $\hat{\beta} = (X'\hat{V}^{-1}X)^{-1}X'\hat{V}^{-1}y$ and $\hat{V}$ is $V$ with $\theta$ replaced by $\hat{\theta}$. Note that the set $\{(i_1,\ldots,i_4):f(i_1,\ldots,i_4) = f_l\}$ does not depend on $\theta$. It follows that $\hat{c}_{j,k}(i_1,\ldots,i_4) = \hat{c}_{j,k,l}$ if $f(i_1,\ldots,i_4) = f_l$, $1 \leq l \leq L$, where $\hat{c}_{j,k,l}$ is given by (14) with $B$ replaced by $\hat{B}$. Here $\hat{B}_{j,i_1,i_2}$ is $B_{j,i_1,i_2}$ with $\theta$ replaced by $\hat{\theta}$, and so forth. This is the observed part.

On the other hand, $\mathcal{I}_{1,2,jk}$ depends only on $\theta$ and, therefore, can be estimated by replacing $\theta$ by $\hat{\theta}$. The result, denoted by $\hat{\mathcal{I}}_{1,2,jk}$, is the estimated part.

An estimator of $\mathcal{I}_{1,jk}$ is then $\hat{\mathcal{I}}_{1,1,jk} + \hat{\mathcal{I}}_{1,2,jk}$; hence an estimator of $\mathcal{I}_1$ is given by $\hat{\mathcal{I}}_1 = \hat{\mathcal{I}}_{1,1} + \hat{\mathcal{I}}_{1,2}$, where $\hat{\mathcal{I}}_{1,r} = (\hat{\mathcal{I}}_{1,r,jk})_{0 \le j,k \le s}$, $r = 1, 2$. Because the estimator consists partially of an observed form and partially of an estimated one, it is called a partially observed quasi-information matrix (POQUIM).

This is exactly where the decomposition (4) came from. We now use another simple example to illustrate the POQUIM decomposition, with more examples to come in Section 3.

EXAMPLE 2. Consider a one-way random effects model $y_{ij} = \mu + \alpha_i + \varepsilon_{ij}$, $i = 1, \ldots, m$, $j = 1, \ldots, n$, where $\mu$ is an unknown mean; the random effects $\alpha_1, \ldots, \alpha_m$ are i.i.d. with mean 0 and variance $\sigma_1^2$; the errors $\varepsilon_{ij}$'s are i.i.d. with mean 0 and variance $\sigma_0^2$; and $\alpha$ and $\varepsilon$ are independent. It is, in this case, more convenient to use a double index (i.e., $ij$ instead of $i$). It is easy to show that $f(i_1 j_1, \ldots, i_4 j_4) = 0$ if not $i_1 = \cdots = i_4$; $\kappa_1$ if $i_1 = \cdots = i_4$ but not $j_1 = \cdots = j_4$; and $\kappa_0 + \kappa_1$ if $i_1 = \cdots = i_4$ and $j_1 = \cdots = j_4$. Thus, $L = 2$ [note that $L$ is the number of different functional values of $f(i_1 j_1, \ldots, i_4 j_4)$]. Define the following functions of $\theta$, where $\theta = (\lambda, \gamma_1)'$: $t_0 = 1 - \{\gamma_1/(1 + \gamma_1 n)\} - 1/(1 + \gamma_1 n)mn$, $t_1 = \{(m-1)n/m(1 + \gamma_1 n)\}$ and $t_3 = \{n(1 + \gamma_1 n)^2 - (1 + \gamma_1)^2\}/(n^3 - 1)$. Then the POQUIM is given by $\hat{\mathcal{I}}_{1,kl} = \hat{\mathcal{I}}_{1,1,kl} + \hat{\mathcal{I}}_{1,2,kl}$, $k, l = 0, 1$, where

$$\hat{\mathcal{I}}_{1,1,00} = \frac{\hat{t}_1^2 - \hat{t}_0^2 n}{4\hat{\lambda}^4 n(n^3 - 1)} \left\{ \sum_i \left( \sum_j \hat{u}_{ij} \right)^4 - \sum_{i,j} \hat{u}_{ij}^4 \right\} + \frac{\hat{t}_0^2}{4\hat{\lambda}^4} \sum_{i,j} \hat{u}_{ij}^4,$$

$$\hat{\mathcal{I}}_{1,1,01} = \frac{(m-1)(\hat{t}_1 n - \hat{t}_0)}{4\hat{\lambda}^3(1 + \hat{\gamma}_1 n)^2 m(n^3 - 1)} \left\{ \sum_i \left( \sum_j \hat{u}_{ij} \right)^4 - \sum_{i,j} \hat{u}_{ij}^4 \right\}$$
$$+ \frac{(m-1)\hat{t}_0}{4\hat{\lambda}^3(1 + \hat{\gamma}_1 n)^2 m} \sum_{i,j} \hat{u}_{ij}^4,$$

$$\hat{\mathcal{I}}_{1,1,11} = \frac{(m-1)^2}{4\hat{\lambda}^2(1 + \hat{\gamma}_1 n)^4 m^2} \sum_i \left( \sum_j \hat{u}_{ij} \right)^4;$$

$$\hat{\mathcal{I}}_{1,2,00} = \frac{1}{2\hat{\lambda}^2} \left[ mn - 1 - \frac{3}{2} mn\hat{t}_0^2 \{(1 + \hat{\gamma}_1)^2 - \hat{t}_3\} - \frac{3}{2} m\hat{t}_1^2 \hat{t}_3 \right],$$

$$\hat{\mathcal{I}}_{1,2,01} = \frac{(m-1)n}{2\hat{\lambda}(1 + \hat{\gamma}_1 n)} \left\{ 1 - \left( \frac{3}{2} \right) \frac{(\hat{t}_1 n - \hat{t}_0)\hat{t}_3 + (1 + \hat{\gamma}_1)^2 \hat{t}_0}{1 + \hat{\gamma}_1 n} \right\},$$

$$\hat{\mathcal{I}}_{1,2,11} = -\frac{(m-1)(m-3)n^2}{4m(1+\hat{\gamma}_1 n)^2},$$

$\hat{u}_{ij} = y_{ij} - \bar{y} \cdots$ and the $\hat{t}$'s are the $t$'s with $\theta$ replaced by $\hat{\theta}$, the REML estimator.

COMPUTATIONAL NOTE. The following list outlines a numerical algorithm for POQUIM:

1. Determine the sets of indices $\mathcal{S}_l = \{(i_1, \ldots, i_4) : f(i_1, \ldots, i_4) = f_l\}$, $1 \le l \le L$. Then, for each $(j,k)$, $0 \le j \le k \le s$, do the following.
2. Compute $\hat{c}_{j,k,l}$, $1 \le l \le L$. Note that the denominator in (14) is $|\mathcal{S}_l|$.
3. Compute $\hat{\mathcal{I}}_{1,1,jk} = \sum_{f(i_1, \ldots, i_4) \neq 0} \hat{c}_{j,k}(i_1, \ldots, i_4) \hat{u}_{i_1} \cdots \hat{u}_{i_4}$, where $\hat{c}_{j,k}(i_1, \ldots, i_4)$ is defined the same way as $c_{j,k}(i_1, \ldots, i_4)$ above (14) with $\theta$ replaced by $\hat{\theta}$ and $\hat{u}_i = y_i - x_i'\hat{\beta}$. Note that $\sum_{f(i_1, \ldots, i_4) \neq 0} = \sum_{\mathcal{S}_1} + \cdots + \sum_{\mathcal{S}_L}$.
4. Compute $\hat{\mathcal{I}}_{1,2,jk}$, which is $\mathcal{I}_{1,2,jk}$ with $\theta$ replaced by $\hat{\theta}$. See step 3 for the summation.
5. Let $\hat{\mathcal{I}}_{1,jk} = \hat{\mathcal{I}}_{1,1,jk} + \hat{\mathcal{I}}_{1,2,jk}$.

All except step 1 are fairly straightforward. As for step 1, the sets may be determined as follows. First, the index $(1,1,1,1)$ belongs to $\mathcal{S}_1$. Also compute the vector $v_{1,1,1,1} = (z_{1t} \cdot z_{1t} \cdot z_{1t} \cdot z_{1t})_{0 \le t \le s}$. Then compute the vector $v_{1,1,1,2} = (z_{1t} \cdot z_{1t} \cdot z_{1t} \cdot z_{2t})_{0 \le t \le s}$. If $v_{1,1,1,2} = v_{1,1,1,1}$, the index $(1,1,1,2)$ belongs to $\mathcal{S}_1$; otherwise it belongs to $\mathcal{S}_2$, and so on.

The main theoretical result in this section is the consistency of POQUIM. To state the result, we need some additional notation.

2.2. *Notation.* For a vector $a = (a_l)$, the 1-norm of $a$ is defined by $\|a\|_1 = \sum_l |a_l|$. A sequence of matrices $M$ is bounded from above if $\|M\|$ is bounded; the sequence is bounded from below if $\|M^{-1}\|$ is bounded, where the norm of a matrix $M$ is defined as $\|M\| = \{\lambda_{\max}(M'M)\}^{1/2}$ with $\lambda_{\max}$ representing the largest eigenvalue. A positive definite matrix-valued function $M(\theta)$, which may depend on $N$, is said to be uniformly continuous at $\theta$ if $\|M^{-1/2}(\theta)M(\theta + \Delta)M^{-1/2}(\theta) - I\| \to 0$ as $\Delta \to 0$ uniformly in $N$, where $I$ is an identity matrix of fixed dimension. It is easy to show that $M(\theta)$ is uniformly continuous if and only if for any $\eta > 0$, there is $\delta > 0$ such that $|\Delta| \le \delta$ implies $(1-\eta)M(\theta) \le M(\theta + \Delta) \le (1+\eta)M(\theta)$ for all $N$. An estimator $\hat{M}$ of a positive definite matrix $M$, which may depend on $N$, is consistent if $M^{-1/2}\hat{M}M^{-1/2} - I \to 0$ in probability, where $I$ is an identity matrix of fixed dimension. In the following discussion, all the sequences of numbers (vectors, matrices) depend on $N$, but for notational simplicity the subscript $N$ is suppressed. For example, in condition (iii) of Theorem 2 below $g_0$ means $g_{0,N}$, et cetera.

Recall that $z'_{it}$ is the $i$th row of $Z_t$, $0 \leq t \leq s$, with $Z_0 = I$. Let $w'_i = (z'_{i0}, z'_{i1}, \ldots, z'_{is})$. Define $d_i$ as a vector of the same dimension as $w_i$ such that the $j$th component of $d_i$ is 1 if the corresponding component of $w_i$ is nonzero and the $j$th component of $d_i$ is 0 if the corresponding component of $w_i$ is zero. Note that $d_i$ is an indicator of what random effects and errors are involved in the expression of $y_i$. Let $h_l$ denote the denominator in (14) and let $h_{l_1, l_2}$ be the cardinality of the set of $(i_1, \ldots, i_8)$ such that $f(i_1, \ldots, i_4) = f_{l_1}$, $f(i_5, \ldots, i_8) = f_{l_2}$ and $(d_{i_1} + \cdots + d_{i_4}) \cdot (d_{i_5} + \cdots + d_{i_8}) \neq 0$. Here, recall that for two vectors $a = (a_l)$ and $b = (b_l)$, the dot product is defined as $a \cdot b = \sum_l a_l b_l$.

Also recall that $\mathcal{G} = 2\{\mathrm{tr}(B_j V B_k V)\}_{0 \leq j, k \leq s}$ is the Gaussian information matrix [see the remark below (14)], that is, $\mathcal{I}_1 = \mathcal{G}$ under normality. More generally, let $\tilde{\mathcal{G}}$ denote $\mathcal{G}$ with $\theta$ replaced by $\tilde{\theta}$ as a function of $\tilde{\theta}$. For any $\delta > 0$, define $g_{j,k}(\delta) = \sup_{\tilde{\theta} \in \Theta, |\tilde{\theta} - \theta| \leq \delta} |\tilde{\mathcal{G}}_{j,k} - \mathcal{G}_{j,k}|$, where $M_{j,k}$ denotes the $j, k$ element of a matrix $M$, and define $d_{j,k,l}(\delta) = \sup_{\tilde{\theta} \in \Theta, |\tilde{\theta} - \theta| \leq \delta} |\tilde{c}_{j,k,l} - c_{j,k,l}|$, where $c_{j,k,l}$ is defined by (14) and $\tilde{c}_{j,k,l}$ is $c_{j,k,l}$ with $\theta$ replaced by $\tilde{\theta}$.

Finally, recall that the asymptotic covariance matrix of the REML estimator, $\hat{\theta}$, is given by (2), that is, $\Sigma_\mathrm{R} = \mathcal{I}_2^{-1} \mathcal{I}_1 \mathcal{I}_2^{-1}$, where $\mathcal{I}_1 = \mathrm{Var}(\partial l_\mathrm{R} / \partial \theta)$ is the QUIM defined in Section 1, and $\mathcal{I}_2 = \mathrm{E}(\partial^2 l_\mathrm{R} / \partial \theta \, \partial \theta')$. The POQUIM estimator of $\Sigma_\mathrm{R}$ is defined by $\hat{\Sigma}_\mathrm{R} = \hat{\mathcal{I}}_2^{-1} \hat{\mathcal{I}}_1 \hat{\mathcal{I}}_2^{-1}$, where $\hat{\mathcal{I}}_1$ is the POQUIM and $\hat{\mathcal{I}}_2$ is the estimated $\mathcal{I}_2$ obtained by replacing the variance components in $\mathcal{I}_2$ by their REML estimators.

### 2.3. *Consistency.*

It should be pointed out that the definition of REML estimator in non-Gaussian mixed linear models differs slightly according to several authors. In [22] the REML estimator is defined as the solution to the REML equation; in [12] the REML estimator is defined as the solution to the REML equation plus the requirement that it belong to the parameter space; in [13] the REML estimator is defined as the maximizer of the Gaussian restricted likelihood. In fact, the last showed that, for balanced mixed linear models, such a maximizer is a consistent estimator of $\theta$; for an unbalanced mixed linear model, it showed that a sieved maximizer is consistent. Note that from a practical point of view the sieve puts no restriction on the maximization, because the maximizer is always within a sieve that satisfies the conditions (of Jiang [13], with a suitable constant). Therefore, in the following theorem the REML estimator is understood as the maximizer of the Gaussian restricted likelihood in the sense of Jiang [13] (with the sieves in the unbalanced case; see above). This eliminates any possible confusion as to which solution, or root, to the REML equation to use when there are multiple roots (e.g., [23], Section 8.1).

THEOREM 2. *Suppose that* (i) $\sigma_t^2 > 0$, $0 < \mathrm{var}(\alpha_{t1}^2) < \infty$, $0 \leq t \leq s$; (ii) $|x_i|$, $\|z_{it}\|_1$, $1 \leq t \leq s$, $1 \leq i \leq N$ *are bounded;* (iii) *there is a sequence of diagonal matrices* $G = \mathrm{diag}(g_0, \ldots, g_s)$ *with* $g_j > 0$, $0 \leq j \leq s$, *such that* $G^{-1} \mathcal{G} G^{-1}$

*is bounded from above as well as from below and* $\lambda_{\min}(X'V^{-1}X) \to \infty$; (iv) $(g_j g_k)^{-1} \sum_{l=1}^{L} h_l |c_{j,k,l}|$, $0 \le j, k \le s$, *are bounded and* $(g_j g_k)^{-2} \times \sum_{l_1, l_2=1}^{L} h_{l_1, l_2} |c_{j,k,l_1} c_{j,k,l_2}| \to 0$, $0 \le j, k \le s$; (v) $(g_j g_k)^{-1} g_{j,k}(\delta) \to 0$ *and* $(g_j g_k)^{-1} \sum_{l=1}^{L} h_l d_{j,k,l}(\delta) \to 0$, $0 \le j, k \le s$, *uniformly in* $N$ *as* $\delta \to 0$. *Then the POQUIM* $\hat{\mathcal{I}}_1$ *and the POQUIM estimator* $\hat{\Sigma}_R$ *are both consistent.*

REMARK 1. The first part of condition (iii) (regarding $\mathcal{G}$) is equivalent to the AI[4] condition of Jiang [12, 13], which, together with $\sigma_t^2 > 0$, $0 \le t \le s$, guarantees the consistency of the REML estimator $\hat{\theta}$. Furthermore, condition (iii) ensures the consistency of $\hat{\beta} = (X'\hat{V}^{-1}X)^{-1}X'\hat{V}^{-1}y$, where $\hat{V}$ is $V$ with $\theta$ replaced by $\hat{\theta}$. Finally, by the proof of Lemma 3 in the sequel, it can be shown that the first part of condition (v) [regarding $g_{j,k}(\delta)$] is equivalent to $\tilde{\mathcal{G}}$ being uniformly continuous at $\theta$.

The proof of Theorem 2 is given in Section 8.3.

**3. Examples.** We now consider some examples and show that the conditions of Theorem 2 are satisfied in typical situations of non-Gaussian mixed linear models.

3.1. *A balanced two-way random effects model.* Example 1 was used in Section 1 to illustrate the POQUIM method. We now revisit this example and verify the conditions of Theorem 2.

Condition (i) is satisfied if $\sigma_t^2 > 0$, $t = 0, 1, 2$, and $0 < \text{var}(v_1^2)$, $\text{var}(w_1^2)$, $\text{var}(e_{11}^2) < \infty$.

Condition (ii) is automatically satisfied, because here $x_{ij} = 1$ and $z_{ijt}$, $t = 0, 1, 2$, are vectors with one component equal to 1 and the other components equal to 0. Note that, as in Example 2, it is more convenient to use a double index, $ij$, instead of $i$.

By Jiang [12], condition (iii) is satisfied with $g_0 = \sqrt{mn}$, $g_1 = \sqrt{m}$ and $g_2 = \sqrt{n}$ if $\sigma_t^2 > 0$, $t = 0, 1, 2$, and $m, n \to \infty$. See the remarks below Theorem 2. Note that here $X'V^{-1}X = mn/\lambda(1 + \gamma_1 n + \gamma_2 m)$.

Now consider condition (iv). It is easy to show that $f(i_1 j_1, \ldots, i_4 j_4) = 0$ if not $i_1 = \cdots = i_4$ or $j_1 = \cdots = j_4$; $\kappa_1$ if $i_1 = \cdots = i_4$ but not $j_1 = \cdots = j_4$; $\kappa_2$ if $j_1 = \cdots = j_4$ but not $i_1 = \cdots = i_4$; and $\kappa_0 + \kappa_1 + \kappa_2$ if $i_1 = \cdots = i_4$ and $j_1 = \cdots = j_4$. Thus, $L = 3$. It is easy to show that $h_1 = mn(n^3 - 1)$, $h_2 = nm(m^3 - 1)$, $h_3 = mn$; $|c_{0,0,1}| \propto n^{-3}$, $|c_{0,0,2}| \propto m^{-3}$, $|c_{0,0,3}| \propto 1$; $|c_{0,1,1}| \propto n^{-4}$, $|c_{0,1,2}| \propto m^{-3}n^{-2}$, $|c_{0,1,3}| \propto n^{-2}$; $|c_{0,2,1}| \propto n^{-3}m^{-2}$, $|c_{0,2,2}| \propto m^{-4}$, $|c_{0,2,3}| \propto m^{-2}$; $|c_{1,1,1}| \propto n^{-4}$, $|c_{1,1,2}| \propto m^{-3}n^{-4}$, $|c_{1,1,3}| \propto n^{-4}$; $|c_{1,2,1}| \propto m^{-2}n^{-5}$, $|c_{1,2,2}| \propto n^{-2}m^{-5}$, $|c_{1,2,3}| \propto m^{-2}n^{-2}$; $|c_{2,2,1}| \propto n^{-3}m^{-4}$, $|c_{2,2,2}| \propto m^{-4}$ and $|c_{2,2,3}| \propto m^{-4}$. It follows that the first part of condition (iv) is satisfied as $m, n \to \infty$.

Furthermore, it is easy to show that $h_{1,1} \propto mn^7(m+n)$, $h_{1,2} \propto m^4n^4(m+n)$, $h_{1,3} \propto mn^4(m+n)$, $h_{2,2} \propto m^7n(m+n)$, $h_{2,3} \propto m^4n(m+n)$ and $h_{3,3} \propto mn \times (m+n)$. It follows that $h_{l_1,l_2} \leq c(m^{-1} + n^{-1})h_{l_1}h_{l_2}$, $1 \leq l_1, l_2 \leq 3$. Therefore, we have

$$(g_jg_k)^{-2} \sum_{l_1,l_2=1}^{3} h_{l_1,l_2}|c_{j,k,l_1}c_{j,k,l_2}| \leq c\left(\frac{1}{m} + \frac{1}{n}\right)\left\{(g_jg_k)^{-1}\sum_{l=1}^{3} h_l|c_{j,k,l}|\right\}^2 \longrightarrow 0$$

as $m, n \to \infty$, $0 \leq j, k \leq 2$, using the already verified first part.

As for condition (v), it is easy to show that the derivatives of $c_{j,k,l}$ with respect to $\theta$ are bounded by quantities of the same order as $|c_{j,k,l}|$. It follows that $d_{j,k,l}(\delta)$ is bounded by $\delta$ times a quantity of the same order as $|c_{j,k,l}|$. Thus, the second part of condition (v) is satisfied by the verified first part of condition (iv). By a similar argument, the first part of condition (v) is satisfied.

In conclusion, all the conditions of Theorem 2 are satisfied with $g_0 = \sqrt{mn}$, $g_1 = \sqrt{m}$ and $g_2 = \sqrt{n}$, provided that $\sigma_t^2 > 0$, $t = 0, 1, 2$, $0 < \text{var}(v_1^2)$, $\text{var}(w_1^2)$, $\text{var}(e_{11}^2) < \infty$ and $m, n \to \infty$.

### 3.2. A balanced two-way mixed effects model.

In the previous example the only fixed effect is an unknown mean $\mu$. This time we consider a model that involves more fixed effects. We assume that $y_{ij} = \beta_j + \alpha_i + \varepsilon_{ij}$, $i = 1, \ldots, m$, $j = 1, \ldots, n$, where the $\beta_j$'s are unknown fixed effects, $\alpha_i$'s are i.i.d. random effects with mean 0 and variance $\sigma_1^2$, $\varepsilon_{ij}$'s are i.i.d. errors with mean 0 and variance $\sigma_0^2$, and $\alpha$ and $\varepsilon$ are independent. Again we verify the conditions of Theorem 2.

Condition (i) holds if $\sigma_t^2 > 0$, $t = 0, 1$, and $0 < \text{var}(\alpha_1^2)$, $\text{var}(\varepsilon_{11}^2) < \infty$.

Condition (ii) is automatically satisfied.

By Jiang [12], condition (iii) is satisfied with $g_0 = \sqrt{mn}$ and $g_1 = \sqrt{m}$ as long as $\sigma_t^2 > 0$, $t = 0, 1$, $m \to \infty$ and $n \geq 2$. Note that this result holds regardless of $n \to \infty$ or not.

Now consider (iv). It is easy to show that $f(i_1j_1, \ldots, i_4j_4) = 0$ if not $i_1 = \cdots = i_4$; $\kappa_1$ if $i_1 = \cdots = i_4$ but not $j_1 = \cdots = j_4$; and $\kappa_0 + \kappa_1$ if $i_1 = \cdots = i_4$ and $j_1 = \cdots = j_4$. Thus $L = 2$. It is easy to verify that $h_1 = mn(n^3 - 1)$ and $h_2 = mn$. Furthermore, we have $|c_{0,0,1}| \propto n^{-3}$, $|c_{0,0,2}| \propto 1$, $|c_{0,1,1}| \propto n^{-4}$, $|c_{0,1,2}| \propto n^{-2}$, $|c_{1,1,1}| \propto n^{-4}$ and $|c_{1,1,2}| \propto n^{-4}$. It follows that the first part of condition (iv) is satisfied as $m \to \infty$. Also, we have $h_{1,1} \propto mn^8$, $h_{1,2} \propto mn^5$ and $h_{2,2} \propto mn^2$; hence $h_{l_1,l_2} \leq cm^{-1}h_{l_1}h_{l_2}$, $1 \leq l_1, l_2 \leq 2$. Thus, for the same reason as in the previous subsection, the second part of condition (iv) is satisfied as $m \to \infty$.

By similar arguments as in the previous subsection, condition (v) is satisfied.

In conclusion, all the conditions of Theorem 2 are satisfied with $g_0 = \sqrt{mn}$ and $g_1 = \sqrt{m}$, provided that $\sigma_t^2 > 0$, $t = 1, 2$, $0 < \text{var}(\alpha_1^2)$, $\text{var}(\varepsilon_{11}^2) < \infty$, $m \to \infty$ and $n \geq 2$.

3.3. *An unbalanced nested error regression model.* In the previous examples the data are balanced in the sense that there are equal numbers of observations per cell (e.g., [23], Chapter 4). In this subsection we consider an unbalanced case. The model may be viewed as an extension of Example 2 in Section 2, which can be expressed as $y_{ij} = x'_{ij}\beta + \alpha_i + \varepsilon_{ij}$, $i = 1, \ldots, m$, $j = 1, \ldots, n_i$, where $n_i$ $(n_i \geq 1)$ is the size of the $i$th cluster, $x_{ij}$ is a vector of known covariates, $\beta$ is a $p$-dimensional vector of unknown regression coefficients, $\alpha_i$'s are i.i.d. random effects with mean 0 and variance $\sigma_1^2$, $\varepsilon_{ij}$'s are i.i.d. errors with mean 0 and variance $\sigma_0^2$, and $\alpha$ and $\varepsilon$ are independent. When $x_{ij} = 1$, $\beta = \mu$ and $n_i = n$, $1 \leq i \leq m$, the model reduces to Example 2 of Section 2. Such a model is useful in a number of application areas, including small area estimation (e.g., [3, 8]). Here, once again, we verify the conditions of Theorem 2.

Condition (i) is satisfied provided that $\sigma_r^2 > 0$, $r = 0, 1$, and $0 < \mathrm{var}(\alpha_1^2)$, $\mathrm{var}(\varepsilon_{11}^2) < \infty$. Condition (ii) is satisfied if $|x_{ij}|$ is bounded. By Jiang [12] it can be shown that condition (iii) is satisfied with $g_0 = \sqrt{N}$, where $N = \sum_{i=1}^m n_i$ is the total sample size and $g_1 = \sqrt{m}$, provided that $m \to \infty$, $p$ is bounded, $\limsup(m/N) < 1$ and

$$\liminf\left[\lambda_{\min}\left\{\frac{1}{m}\sum_{i=1}^m\sum_{j=1}^{n_i}(x_{ij} - \bar{x}_{i\cdot})(x_{ij} - \bar{x}_{i\cdot})'\right\} \vee \lambda_{\min}\left\{\frac{1}{m}\sum_{i=1}^m \bar{x}_{i\cdot}\bar{x}'_{i\cdot}\right\}\right] > 0,$$

(15)

where $\bar{x}_{i\cdot} = n_i^{-1}\sum_{j=1}^{n_i} x_{ij}$. Note that $\limsup(m/N) < 1$ ensures that, asymptotically, the random effects and errors can be separated, that is, the variance components are asymptotically identifiable [12]. Although condition (15) can be further weakened, it is more intuitive and satisfied in most cases.

Now consider condition (iv). The function $f(i_1 j_1, \ldots, i_4 j_4)$ has the same expression as in Example 2 of Section 2. Thus we have $L = 2$, $h_1 = \sum_{i=1}^m n_i(n_i^3 - 1)$ and $h_2 = N$. Furthermore, it can be shown that the $i_1 j_1$, $i_2 j_2$ element of $B_0$ is

$$B_{0, i_1 j_1, i_2 j_2} = \frac{1}{2\lambda^2}\left\{\mathbf{1}_{(i_1 = i_2, j_1 = j_2)} - \frac{\gamma_1}{1 + \gamma_1 n_{i_1}}\mathbf{1}_{(i_1 = i_2)}\right.$$
$$\left. - \left(x_{i_1 j_1} - \frac{\gamma_1 n_{i_1}}{1 + \gamma_1 n_{i_1}}\bar{x}_{i_1\cdot}\right)' D^{-1}\left(x_{i_2 j_2} - \frac{\gamma_1 n_{i_2}}{1 + \gamma_1 n_{i_2}}\bar{x}_{i_2\cdot}\right)\right\},$$

where $D = \sum_{i=1}^m X'_i D_i X_i$ with $X_i = (x'_{ij})_{1 \leq j \leq n_i}$ and $D_i = I_{n_i} - \gamma_1(1 + \gamma_1 n_i)^{-1} J_{n_i}$. Similarly, the $i_1 j_1$, $i_2 j_2$ element of $B_1$ is

$$B_{1, i_1 j_1, i_2 j_2} = \frac{1}{2\lambda}\left\{\frac{\mathbf{1}_{(i_1 = i_2)}}{(1 + \gamma_1 n_{i_1})(1 + \gamma_1 n_{i_2})}\right.$$
$$\left. - \frac{n_{i_1}}{(1 + \gamma_1 n_{i_1})^2}\left(x_{i_2 j_2} - \frac{\gamma_1 n_{i_2}}{1 + \gamma_1 n_{i_2}}\bar{x}_{i_2\cdot}\right)' D^{-1}\bar{x}_{i_1\cdot}.\right.$$

$$-\frac{n_{i_2}}{(1+\gamma_1 n_{i_2})^2}\Big(x_{i_1 j_1}-\frac{\gamma_1 n_{i_1}}{1+\gamma_1 n_{i_1}}\bar{x}_{i_1\cdot}\Big)' D^{-1}\bar{x}_{i_2\cdot}$$

$$+\Big(x_{i_1 j_1}-\frac{\gamma_1 n_{i_1}}{1+\gamma_1 n_{i_1}}\bar{x}_{i_1\cdot}\Big)' D^{-1}Q D^{-1}\Big(x_{i_2 j_2}-\frac{\gamma_1 n_{i_2}}{1+\gamma_1 n_{i_2}}\bar{x}_{i_2\cdot}\Big)\Big\},$$

where $Q=\sum_{i=1}^{m}\{n_i/(1+\gamma_1 n_i)\}^2\bar{x}_{i\cdot}\bar{x}_{i\cdot}'$. Thus, it can be shown that $|c_{0,0,1}|\propto N/\sum_{i=1}^{m}n_i(n_i^3-1)$, $|c_{0,0,2}|\propto 1$, $|c_{0,1,1}|,|c_{1,1,1}|\propto m/\sum_{i=1}^{m}n_i(n_i^3-1)$ and $|c_{0,1,2}|,|c_{1,1,2}|\propto m/N$. It follows that the first part of condition (iv) is satisfied. Furthermore, it can be shown that $|h_{l_1,l_2}|\le ch_{l_1}h_{l_2}\sum_{i=1}^{m}n_i^{a_1+a_2}/(\sum_{i=1}^{m}n_i^{a_1})\times(\sum_{i=1}^{m}n_i^{a_2})$, $l_1,l_2=1,2$, provided that $\limsup(m/N)<1$, where $c$ is a constant and $a_r=(3-l_r)^2$, $r=1,2$. Note that here we use the fact that, by Hölder's inequality and the fact that $n_i\ge 1$, it can be shown that $(\sum_{i=1}^{m}n_i)/(\sum_{i=1}^{m}n_i^4)\le (m/N)^{3/4}$, which implies that $h_1\propto\sum_{i=1}^{m}n_i^4$, because $\limsup(m/N)<1$. Thus, by the first part of condition (iv), the second part of condition (iv) is satisfied, provided that

$$(16)\qquad\frac{\sum_{i=1}^{m}n_i^{a+b}}{(\sum_{i=1}^{m}n_i^a)(\sum_{i=1}^{m}n_i^b)}\longrightarrow 0,\qquad a,b=1\text{ or }4.$$

Note that, for example, if the $n_i$'s are bounded, then the left-hand side of (16) is $O(m^{-1})$.

Finally, condition (v) can be verified using the same arguments as in part (v) of the previous two subsections.

In conclusion, all the conditions of Theorem 2 are satisfied with $g_0=\sqrt{N}$ and $g_1=\sqrt{m}$, provided that $\sigma_r^2>0$, $r=0,1$, $0<\mathrm{var}(\alpha_1^2)$, $\mathrm{var}(\varepsilon_{11}^2)<\infty$, $m\to\infty$, $p$ is bounded, $\limsup(m/N)<1$, and (15) and (16) hold. Note that the conditions *do not* include that $n_i\to\infty$. In fact, in most practical situations the $n_i$'s are small (e.g., Ghosh and Rao [8]).

### 3.4. *A random intercept/slope model.*

So far in the examples the number of different functional values for $f(i_1,\dots,i_4)$, defined by (8), is bounded, that is, $L$ is bounded. We now consider a case in which $L$ increases with $N$.

Suppose that two measures are collected from each of $m$ patients, once before and once after a surgery, but, because of the availability of patients, the measures are made at different times after the time of the surgery. It is thought that the recovery is a linear function of time, but the slope depends on the individual patient. For the $i$th patient, the baseline measure made before the surgery can be expressed as $\beta_0+a_i$, where $\beta_0$ is an unknown mean and $a_i$ is a random effect; after the surgery, a measure is collected at time $t_i$ from the surgery and the improvement can be expressed as $(\beta_1+b_i)t_i$ on top of the baseline, where $\beta_1$ is an unknown parameter and $b_i$ is another random effect. Of course, each time there is a random measurement error. This model can be expressed as $y_i=\beta_0+a_i+e_i$, $y_{m+i}=\beta_0+\beta_1 t_i+a_i+b_i t_i+e_{m+i}$, $i = 1, \ldots, m$, where $y_i$ and $y_{m+i}$ correspond to the measurements from the $i$th patient before and after the surgery. It is assumed that the $a_i$'s are i.i.d. with mean 0 and variance $\sigma_1^2$, the $b_i$'s are i.i.d. with mean 0 and variance $\sigma_2^2$, the $e_i$'s are i.i.d. with mean 0 and variance $\sigma_0^2$, and $a$, $b$, $e$ are independent (see the discussion in Section 7).

Now consider the conditions of Theorem 2. Condition (i) is satisfied if $\sigma_t^2 > 0$, $t = 0, 1, 2$, and $0 < \mathrm{var}(a_1^2)$, $\mathrm{var}(b_1^2)$, $\mathrm{var}(e_1^2) < \infty$. Condition (ii) is satisfied provided that the $t_i$'s are bounded. By Jiang [12], it can be shown that condition (iii) is satisfied with $g_j = \sqrt{m}$, $j = 0, 1, 2$, provided that $m \to \infty$ and the $t_i$'s are bounded from above and away from zero.

Now consider condition (iv). For $1 \le i \le N = 2m$, write $i = 2(l-1) + r$, where $1 \le l \le m$ and $r = 1, 2$. Then it can be shown that the $i_1, i_2$ element of $B_j$ can be expressed as $O(1)\mathbf{1}_{(l_1 = l_2)} + O(m^{-1})$, $j = 0, 1, 2$. For simplicity, assume that the $t_i$'s are all different. Then it is easy to see that $f(i_1, \ldots, i_4) = 0$ if not $l_1 = \cdots = l_4$; $\kappa_1$ if $l_1 = \cdots = l_4$ but not $r_1 = \cdots = r_4$; $\kappa_0 + \kappa_1$ if $l_1 = \cdots = l_4$ and $r_1 = \cdots = r_4 = 1$; and $\kappa_0 + \kappa_1 + t_l^2\kappa_2$ if $l_1 = \cdots = l_4 = l$ and $r_1 = \cdots = r_4 = 2$, $1 \le l \le m$. Thus, we have, in particular, $L = m + 2$. It is easy to show that $h_1 = 14m$, $h_2 = m$ and $h_{2+l} = 1$, $1 \le l \le m$. Furthermore, we have $|c_{j,k,1}| = O(1)$, $1 \le l \le m + 2$. It follows that $(g_j g_k)^{-1}\{h_1|c_{j,k,1}| + h_2|c_{j,k,2}| + \sum_{l=1}^{m} h_{2+l}|c_{j,k,2+l}|\} = m^{-1}O(m) = O(1)$; hence the first part of condition (iv) is satisfied. Similarly, we have $h_{1,1} = O(m)$, $h_{1,2} = O(m)$, $h_{1,2+l} = O(1)$, $1 \le l \le m$, $h_{2,2} = m$, $h_{2,2+l} = 1$, $1 \le l \le m$, and $h_{2+l,2+l'} = \mathbf{1}_{(l=l')}$, $1 \le l, l' \le m$. Thus, we have $(g_j g_k)^{-2} \sum_{l_1,l_2=1}^{L} h_{l_1,l_2}|c_{j,k,l_1}c_{j,k,l_2}| = O(m^{-1})$; hence the second part of condition (iv) is satisfied. By similar arguments as in the previous examples, condition (v) is satisfied.

In conclusion, all the conditions of Theorem 2 are satisfied with $g_j = \sqrt{m}$, $j = 0, 1, 2$, provided that $\sigma_t^2 > 0$, $t = 0, 1, 2$, $0 < \mathrm{var}(a_1^2)$, $\mathrm{var}(b_1^2)$, $\mathrm{var}(e_1^2) < \infty$, $m \to \infty$, and the $t_i$'s are bounded from above and away from zero, and are different. The last condition that the $t_i$'s are all different is only for technical convenience (otherwise $L$ may be less than $m + 2$, but the conditions can be verified similarly).

**4. POQUIM for ML.** In this section we derive POQUIM for ML estimation. Under model (1) and normality, the log-likelihood for estimating $\beta$ and $\theta$ is given by

$$(17) \qquad l(\beta, \theta) = c - \tfrac{1}{2}\{\log(|V|) + (y - X\beta)'V^{-1}(y - X\beta)\},$$

where $c$ is a constant. If normality does not hold, (17) is considered the quasi-log-likelihood. It is easy to show that $\partial l/\partial \beta = X'V^{-1}u$, and $\partial l/\partial \theta_j = u'C_j u - c_j$, $0 \le j \le s$, where $C_0 = (2\lambda)^{-1}V^{-1}$, $C_j = (\lambda/2)V^{-1}Z_jZ_j'V^{-1}$, $1 \le j \le 1$, $c_0 = N/2\lambda$, $c_j = (\lambda/2)\mathrm{tr}(V^{-1}Z_jZ_j')$, $1 \le j \le s$, and again, $u = y - X\beta$. Note that $c_j = \mathrm{E}(u'C_j u)$, $0 \le j \le s$. Let $V^{-1}X_j = q_j = (q_{j,i})_{1 \le i \le N}$, where $X_j$ is the $j$th column of $X$.

Using the expression (17), it can be shown that

$$\text{cov}\left(\frac{\partial l}{\partial \beta_j}, \frac{\partial l}{\partial \beta_k}\right) = X_j' V^{-1} X_k, \qquad 1 \le j, k \le p.$$

Next, similar to Lemma 1, the following equations can be easily derived.

LEMMA 2.  *We have*

$$\text{cov}(u_{i_1}, u_{i_2} u_{i_3}) = \sum_{t=0}^{s} \text{E}(\alpha_{t1}^3) z_{i_1 t} \cdot z_{i_2 t} \cdot z_{i_3 t}, \tag{18}$$

$$\text{cov}\left(\frac{\partial l}{\partial \beta_j}, \frac{\partial l}{\partial \theta_k}\right) = \sum_{t=0}^{s} \text{E}(\alpha_{t1}^3) \sum_{l=1}^{m_t} (X_j' V^{-1} z_{tl})(z_{tl}' C_k z_{tl}).$$

Write $t(i_1, i_2, i_3) = \text{E}(u_{i_1} u_{i_2} u_{i_3})$, which is the right-hand side of (18). Let $t_l$, $1 \le l \le K$, be the different functional values of $t(i_1, i_2, i_3)$ [as functions of the third moments; see (18)]. Then, by similar arguments as in the previous section, it can be shown that

$$\text{cov}\left(\frac{\partial l}{\partial \beta_j}, \frac{\partial l}{\partial \theta_k}\right) = \text{E}\left\{\sum_{t(i_1,i_2,i_3) \ne 0} c_{1,j,k}(i_1, i_2, i_3) u_{i_1} u_{i_2} u_{i_3}\right\},$$

where $c_{1,j,k}(i_1, i_2, i_3) = c_{1,j,k,l}$ if $t(i_1, i_2, i_3) = t_l$, $1 \le l \le K$, with

$$c_{1,j,k,l} = \frac{1}{|\{t(i_1,i_2,i_3) = t_l\}|} \sum_{t(i_1,i_2,i_3)=t_l} q_{j,i_1} C_{k,i_2,i_3}. \tag{19}$$

Furthermore, recall that $f(i_1, \ldots, i_4)$ is defined by (8). Then, similar to the previous section, define $c_{2,j,k}(i_1, \ldots, i_4) = c_{2,j,k,l}$ if $f(i_1, \ldots, i_4) = f_l$, $1 \le l \le L$, with

$$c_{2,j,k,l} = \frac{1}{|\{f(i_1,\ldots,i_4) = f_l\}|} \sum_{f(i_1,\ldots,i_4)=f_l} C_{j,i_1,i_2} C_{k,i_3,i_4}. \tag{20}$$

Then we have similar expressions for $\text{cov}(\partial l/\partial \theta_j, \partial l/\partial \theta_k)$ (with the only difference from Section 2.1 being that $B$ is replaced by $C$). We summarize the results as follows. As before, write $\psi = (\beta' \theta')'$ and, again, write the QUIM as $\mathcal{I}_1 = \text{Var}(\partial l/\partial \psi) = (\mathcal{I}_{1,jk})_{1 \le j,k \le p+s+1}$.

THEOREM 3.  *For any non-Gaussian mixed linear model* (1), *we have*

$$\mathcal{I}_{1,jk} = X_j' V^{-1} X_k = \mathcal{I}_{1,2,jk}, \tag{21}$$

*that is, $\mathcal{I}_{1,1,jk} = 0$, $1 \le j, k \le p$;*

$$\mathcal{I}_{1,j(p+k+1)} = \sum_{t=0}^{s} \mathrm{E}(\alpha_{t1}^3) \sum_{l=1}^{m_t} (X_j' V^{-1} z_{tl})(z_{tl}' C_k z_{tl})$$

$$(22) \qquad = \mathrm{E}\left\{ \sum_{t(i_1,i_2,i_3) \ne 0} c_{1,j,k}(i_1, i_2, i_3) u_{i_1} u_{i_2} u_{i_3} \right\}$$

$$= \mathcal{I}_{1,1,j(p+k+1)},$$

*that is, $\mathcal{I}_{1,2,j(p+k+1)} = 0$, $1 \le j \le p$, $0 \le k \le s$; and*

$$\mathcal{I}_{1,(p+j+1)(p+k+1)} = 2\,\mathrm{tr}(C_j V C_k V) + \sum_{t=0}^{s} \kappa_t \sum_{l=1}^{m_t} (z_{tl}' C_j z_{tl})(z_{tl}' C_k z_{tl})$$

$$= \mathrm{E}\left\{ \sum_{f(i_1,\ldots,i_4) \ne 0} c_{2,j,k}(i_1,\ldots,i_4) u_{i_1} \cdots u_{i_4} \right\}$$

$$(23) \qquad + \left\{ 2\,\mathrm{tr}(C_j V C_k V) \right.$$

$$\left. - 3\lambda^2 \sum_{f(i_1,\ldots,i_4) \ne 0} c_{2,j,k}(i_1,\ldots,i_4) \Gamma(i_1,i_3) \Gamma(i_2,i_4) \right\}$$

$$= \mathcal{I}_{1,1,(p+j+1)(p+k+1)} + \mathcal{I}_{1,2,(p+j+1)(p+k+1)},$$

$0 \le j, k \le s$, *where $c_{1,j,k}(i_1,i_2,i_3)$ and $c_{2,j,k}(i_1,\ldots,i_4)$ are defined by* (19) *and* (20), *respectively.*

Similar to Section 2, the POQUIM is given by $\hat{\mathcal{I}}_1 = (\hat{\mathcal{I}}_{1,jk})_{1 \le j,k \le p+s+1}$, where $\hat{\mathcal{I}}_{1,jk} = \hat{\mathcal{I}}_{1,1,jk} + \hat{\mathcal{I}}_{1,2,jk}$, $\hat{\mathcal{I}}_{1,1,jk}$ is the observed part obtained by taking the expectation sign out of $\mathcal{I}_{1,1,jk}$, if there is one, and then replacing the parameters involved by their ML estimators; and $\hat{\mathcal{I}}_{1,2,jk}$ is the estimated part obtained by replacing $\theta$ by $\hat{\theta}$, the ML estimator, in $\mathcal{I}_{1,2,jk}$, if the latter is nonzero. Let $\hat{\psi} = (\hat{\beta}' \hat{\theta}')'$ be the ML estimator of $\psi$. Then, according to the discussion in Section 1, the ACM of $\hat{\psi}$ is given by $\Sigma = \mathcal{I}_2^{-1} \mathcal{I}_1 \mathcal{I}_2^{-1}$, where $\mathcal{I}_2 = \mathrm{E}(\partial^2 l / \partial \psi \, \partial \psi')$. Thus, the POQUIM estimator of $\Sigma$ is given by $\hat{\Sigma} = \hat{\mathcal{I}}_2^{-1} \hat{\mathcal{I}}_1 \hat{\mathcal{I}}_2^{-1}$, where $\hat{\mathcal{I}}_1$ is the POQUIM and $\hat{\mathcal{I}}_2$ is $\mathcal{I}_2$ with $\theta$ replaced by $\hat{\theta}$. Similar to Theorem 2, sufficient conditions can be given for the consistency of $\hat{\mathcal{I}}_1$ and $\hat{\Sigma}$. The details are omitted.

We now use a simple example to illustrate the POQUIM for ML given by (21)–(23).

EXAMPLE 2 (continued).    Here $p = s = 1$. It is easy to show that $\hat{\mathcal{I}}_{1,11} = mn/\hat{\lambda}(1 + \hat{\gamma}_1 n)$,

$$\hat{\mathcal{I}}_{1,12} = \frac{1}{2\hat{\lambda}^3(1 + \hat{\gamma}_1 n)^2}\left[\frac{1 - \hat{\gamma}_1}{n + 1}\sum_i\left(\sum_j \hat{u}_{ij}\right)^3 + \left\{\hat{\gamma}_1 n + \frac{(1 - \hat{\gamma}_1)n}{n + 1}\right\}\sum_{i,j}\hat{u}_{ij}^3\right]$$

and $\hat{\mathcal{I}}_{1,13} = \{1/2\hat{\lambda}^2(1 + \hat{\gamma}_1 n)^3\}\sum_i(\sum_j \hat{u}_{ij})^3$, where $\hat{\lambda}$, $\hat{\gamma}_1$ are the ML estimators. Furthermore, we have $\hat{\mathcal{I}}_{1,(j+2)(k+2)} = \hat{\mathcal{I}}_{1,1,(j+2)(k+2)} + \hat{\mathcal{I}}_{1,2,(j+2)(k+2)}$, $j, k = 0, 1,$ where

$$\hat{\mathcal{I}}_{1,1,22} = \frac{n\{n - (\hat{\gamma}_1 n + 1 - \hat{\gamma}_1)^2\}}{4\hat{\lambda}^4(1 + \hat{\gamma}_1 n)^2 n(n^3 - 1)}\left\{\sum_i\left(\sum_j \hat{u}_{ij}\right)^4 - \sum_{i,j}\hat{u}_{ij}^4\right\}$$
$$+ \frac{(\hat{\gamma}_1 n + 1 - \hat{\gamma}_1)^2}{4\hat{\lambda}^4(1 + \hat{\gamma}_1 n)^2}\sum_{i,j}\hat{u}_{ij}^4,$$

$$\hat{\mathcal{I}}_{1,1,23} = \frac{n + 1 - \hat{\gamma}_1}{4\hat{\lambda}^3(1 + \hat{\gamma}_1 n)^3(n^2 + n + 1)}\left\{\sum_i\left(\sum_j \hat{u}_{ij}\right)^4 - \sum_{i,j}\hat{u}_{ij}^4\right\}$$
$$+ \frac{\hat{\gamma}_1 n + 1 - \hat{\gamma}_1}{4\hat{\lambda}^3(1 + \hat{\gamma}_1 n)^3}\sum_{i,j}\hat{u}_{ij}^4,$$

$$\hat{\mathcal{I}}_{1,2,22} = \frac{mn}{2\hat{\lambda}^2}\left[1 + \left(\frac{3}{2}\right)\frac{n(1 + \hat{\gamma}_1 n)^2 - (1 + \hat{\gamma}_1)^2}{n^3 - 1}\right.$$
$$\times\left\{\left(\frac{\hat{\gamma}_1 n + 1 - \hat{\gamma}_1}{1 + \hat{\gamma}_1 n}\right)^2 - \frac{n}{(1 + \hat{\gamma}_1 n)^2}\right\}$$
$$\left.- \left(\frac{3}{2}\right)\left(\frac{\hat{\gamma}_1 n + 1 - \hat{\gamma}_1}{1 + \hat{\gamma}_1 n}\right)^2(1 + \hat{\gamma}_1)^2\right],$$

$$\hat{\mathcal{I}}_{1,2,23} = \frac{mn}{2\hat{\lambda}(1 + \hat{\gamma}_1 n)}\left[1 - \left(\frac{3}{2}\right)\frac{(n + 1 - \hat{\gamma}_1)\{n(1 + \hat{\gamma}_1 n)^2 - (1 + \hat{\gamma}_1)^2\}}{(1 + \hat{\gamma}_1 n)^2(n^2 + n + 1)}\right.$$
$$\left.- \left(\frac{3}{2}\right)\frac{(\hat{\gamma}_1 n + 1 - \hat{\gamma}_1)(1 + \hat{\gamma}_1)^2}{(1 + \hat{\gamma}_1 n)^2}\right],$$

$$\hat{\mathcal{I}}_{1,1,33} = \left\{\frac{1}{4\hat{\lambda}^2(1 + \hat{\gamma}_1 n)^4}\right\}\sum_i\left(\sum_j \hat{u}_{ij}\right)^4 \quad\text{and}\quad \hat{\mathcal{I}}_{1,2,33} = \frac{mn^2}{4(1 + \hat{\gamma}_1 n)^2}.$$

**5. Robust dispersion tests.**    In this section we consider an application of the results on POQUIM to robust dispersion tests in mixed linear models. The tests considered here are robust in the sense that they do not require normality. A dispersion test may be regarding both the fixed effects and the variance components or only the variance components, and both ML and REML estimators may be used in such a test. To be more specific, here we consider dispersion tests regarding only the variance components based on the REML estimators.

Consider the following general hypothesis regarding $\theta$ in model (1):

$$\text{H}_0 : K'\theta = \varphi, \tag{24}$$

where $\varphi$ is a specified vector and $K$ is a known $(s+1) \times r$ matrix with rank$(K) = r$. We assume that the REML estimator $\hat{\theta}$ is asymptotically normal with mean 0 and ACM $\Sigma_{\text{R}}$, that is,

$$\Sigma_{\text{R}}^{-1/2}(\hat{\theta} - \theta) \longrightarrow N(0, I_{s+1}) \qquad \text{in distribution.} \tag{25}$$

Sufficient conditions for (25) can be found in, for example, [12]. It is then easy to show that, under the null hypothesis (24),

$$(K'\hat{\theta} - \varphi)'(K'\Sigma_{\text{R}}K)^{-1}(K'\hat{\theta} - \varphi) \longrightarrow \chi_r^2 \qquad \text{in distribution.} \tag{26}$$

We then replace $\Sigma_{\text{R}}$ by its POQUIM estimator $\hat{\Sigma}_{\text{R}}$ of Section 2 to obtain the test statistic

$$\hat{\chi}^2 = (K'\hat{\theta} - \varphi)'(K'\hat{\Sigma}_{\text{R}}K)^{-1}(K'\hat{\theta} - \varphi). \tag{27}$$

The following theorem states that $\hat{\chi}^2$ has the same asymptotic null distribution as (26).

THEOREM 4. *Suppose that the conditions of Theorem 2 are satisfied. Furthermore suppose that (25) holds. Then, under the null hypothesis, $\hat{\chi}^2 \rightarrow \chi_r^2$ in distribution.*

In cases where some components of $\theta$ are specified under the null hypothesis, it is customary to use these specified values, instead of the estimators, in the POQUIM estimator. Under the null hypothesis this may improve the accuracy of the POQUIM estimator, although the difference is expected to be small in large samples (because of the consistency of $\hat{\theta}$; e.g., [12]). It is easy to see, by examing the proofs of Theorems 2 and 4, that the same conclusion of Theorem 4 holds after such a modification. (Note that the only property of $\hat{\theta}$ used in the proof of Theorem 2 is its consistency.) We consider a simple example.

EXAMPLE 2 (continued). Suppose that one wishes to test the hypothesis $\text{H}_0 : \gamma_1 = 1$, that is, the variance contribution due to the random effects is the same as that due to the errors. Note that in this case $\theta = (\lambda, \gamma_1)'$, so the null hypothesis corresponds to (24) with $K = (0, 1)'$ and $\varphi = 1$. Furthermore, we have $K'\Sigma_{\text{R}}K = \Sigma_{\text{R},11}$, which is the asymptotic variance of $\hat{\gamma}_1$, the REML estimator of $\gamma_1$. Thus, the test statistic is $\hat{\chi}^2 = (\hat{\gamma}_1 - 1)^2 / \hat{\Sigma}_{R,11}$, where $\hat{\Sigma}_{R,11}$ is the POQUIM estimator of $\Sigma_{R,11}$ (see Section 2). It is easy to show that

$$(28) \qquad \hat{\Sigma}_{R,11} = \frac{\hat{\mathcal{I}}_{1,11}\hat{\mathcal{I}}_{2,00}^2 - 2\hat{\mathcal{I}}_{1,01}\hat{\mathcal{I}}_{2,00}\hat{\mathcal{I}}_{2,01} + \hat{\mathcal{I}}_{1,00}\hat{\mathcal{I}}_{2,01}^2}{(\hat{\mathcal{I}}_{2,00}\hat{\mathcal{I}}_{2,11} - \hat{\mathcal{I}}_{2,01}^2)^2},$$

where $\hat{\mathcal{I}}_{1,jk} = \hat{\mathcal{I}}_{1,1,jk} + \hat{\mathcal{I}}_{1,2,jk}$, $j,k = 0,1$, and $\hat{\mathcal{I}}_{1,r,jk}$, $r = 1,2$, are given in Example 2 in Section 2, but with $\hat{\gamma}_1$ replaced by 1, its value under $H_0$; furthermore, we have $\hat{\mathcal{I}}_{2,00} = -(mn-1)/2\hat{\lambda}^2$, $\hat{\mathcal{I}}_{2,01} = -(m-1)n/2\hat{\lambda}(1 + \hat{\gamma}_1 n)$ and $\hat{\mathcal{I}}_{2,11} = -(m-1)n^2/2(1 + \hat{\gamma}_1 n)^2$, again with $\hat{\gamma}_1$ replaced by 1, where $\hat{\lambda}$ is the REML estimator of $\lambda$. The asymptotic null distribution is $\chi_1^2$. In the next section the finite sample performance of this test will be investigated.

**6. Simulations.** In this section we consider two simulated examples. The goal is to study the finite sample performance of POQUIM, whose large sample properties were studied in Section 2 (Theorem 2) and later in Section 4 in the context of the robust dispersion test (Theorem 4). The latter will be the focus of our simulation study.

The first example is the one-way random effects model considered in Example 2. Note that this model is a special case of the unbalanced nested error regression model of Section 3.3. However, by restricting to the balanced case we are able to make a direct comparison with the delete-group jackknife method [1, 2].

The second example is a balanced two-way random effects model. Note that the jackknife method does not apply to this case. In fact, when the random effects and errors are not normal or symmetric, POQUIM is the only method that is known to apply, by Section 2, at least in large samples. Now our goal is to investigate its finite sample performance.

6.1. *A balanced one-way random effects model.* Consider once again Example 2 (continued) in Section 5, where the hypothesis to be tested is $H_0$: $\gamma_1 = 1$. For example, such a test may be of genetic interest, which corresponds to $H_0$: $h^2 = 2$, where $h^2 = 4\sigma_1^2/(\sigma_0^2 + \sigma_1^2)$ is the heritability. We consider a test based on REML estimation of the variance components. More specifically, we are interested in the situation when $m$ is increasing while $n$ remains fixed. Therefore, the following sample size configurations are considered: Case I, $m = 50$, $n = 2$; Case II, $m = 400$, $n = 2$. Case I represents a moderate sample size, while Case II represents a large sample size. In addition, we would like to investigate different cases in which normality and symmetry may or may not hold. Therefore, the following combinations of distributions for the random effects and errors are considered: Case i, Normal–Normal; Case ii, DE–NM($-2, 2, 0.5$), where DE represents the double exponential distribution and NM($\mu_1, \mu_2, \rho$) denotes the mixture of two normal distributions with means $\mu_1$, $\mu_2$, variance 1 and mixing probability $\rho$ [i.e., the probabilities $1 - \rho$ and $\rho$ correspond to $N(\mu_1, 1)$ and $N(\mu_2, 1)$, resp.]; and Case iii, CE–NM$(-4, 1, 0.2)$, where CE represents the centralized exponential distribution, that is, the distribution of $X - 1$, where $X \sim$ Exponential(1). Note that in Case ii the distributions are not normal but symmetric, while in Case iii the distributions are not even symmetric—a further departure from normality. Also note that all these distributions have mean 0. They are standardized so that the distributions of the random effects and errors have variances $\sigma_1^2$ and $\sigma_0^2$, respectively. The true value of $\mu$ is set to 1.0. The true value of $\sigma_0^2$ is also chosen as 1.0.

According to Section 5, the $\chi^2$-test statistic is given by

$$\chi^2 = \frac{(\hat{\gamma}_1 - 1)^2}{\hat{\Sigma}_{R,11}}, \tag{29}$$

where $\hat{\gamma}_1$ is the REML estimator of $\gamma_1$, $\hat{\Sigma}_{R,11}$ is given by (28) and

$$\hat{\lambda} = \frac{1}{mn - 1} \left( \text{SSE} + \frac{\text{SSA}}{n + 1} \right). \tag{30}$$

Here

$$\text{SSE} = \sum_{i=1}^{m} \sum_{j=1}^{n} (y_{ij} - \bar{y}_{i\cdot})^2$$

and

$$\text{SSA} = n \sum_{i=1}^{m} (\bar{y}_{i\cdot} - \bar{y}_{\cdot\cdot})^2,$$

with $\bar{y}_{i\cdot} = n^{-1} \sum_{j=1}^{n} y_{ij}$ and $\bar{y}_{\cdot\cdot} = (mn)^{-1} \sum_{i=1}^{m} \sum_{j=1}^{n} y_{ij}$. Note that (30) is the REML estimator of $\lambda$ under the null.

Arvesen [1] proposed a delete-group jackknife method and established consistency, using $U$-statistics. Furthermore, Arvesen and Schmitz [2] provided simulation results. The delete-group jackknife applies to cases where data can be divided into i.i.d. groups, such as the current situation. We refer to this method as jackknife. The method is briefly described as follows. Let $X_1, \ldots, X_m$ be i.i.d. observations and let $\theta$ be an unknown parameter. Let $\hat{\theta}$ be an estimator of $\theta$ based on all the observations and let $\hat{\theta}_{-i}$ be the estimator based on all but the $i$th observation, obtained otherwise the same way as $\hat{\theta}$. Define $\hat{\theta}_i = m\hat{\theta} - (m-1)\hat{\theta}_{-i}$, $1 \leq i \leq m$. The jackknife estimator of $\theta$ is defined as $\hat{\theta}_{\text{jack}} = m^{-1} \sum_{i=1}^{m} \hat{\theta}_i$, that is, the average of the $\hat{\theta}_i$'s.

Now consider the one-way random effects model of Example 2. Instead of deleting the $i$th observation, one deletes the $i$th group consisting of the observations $y_{ij}$, $j = 1, \ldots, n$. A dispersion test that is often of genetic interest is $\text{H}_0$: $\gamma_1 = \gamma_{10}$, which corresponds to $\text{H}_0$: $h^2 = 4\gamma_{10}/(1 + \gamma_{10})$, where $h^2 = 4\sigma_1^2/(\sigma_0^2 + \sigma_1^2)$ is the heritability. Arvesen and Schmitz proposed use of the jackknife estimator with a transformation. Let $\theta = \log(1 + \gamma_1 n)$ and $\hat{\theta} = \log(\text{MSA}/\text{MSE})$, where MSA and MSE are the between and within group mean squares. A test of $H_0$ will be based on

$$(31) \qquad t = \frac{\sqrt{m}(\hat{\theta}_{\text{jack}} - \theta)}{\sqrt{(m-1)^{-1}\sum_{i=1}^m(\hat{\theta}_i - \hat{\theta}_{\text{jack}})^2}},$$

which is expected to have an asymptotic $t_{m-1}$ null distribution [1].

To make a fair comparison, we note that a test of $\chi^2$ type is *omnibus* rather than *directional* (e.g., [21], Section 1.1). In other words, a $\chi^2$ test is typically used in situations of two-sided hypotheses. On the other hand, a $t$-test is appropriate to both one- and two-sided hypotheses. Therefore, we consider testing $H_0$: $\gamma_1 = 1$ against $H_1$: $\gamma_1 \neq 1$. For each simulated data set, the test statistics (29) and (31) are computed. The simulated sizes that correspond to the usual nominal levels 0.01, 0.05 and 0.10 are reported in Table 1. Furthermore, the simulated powers at a number of alternatives, namely, $\gamma_1 = 0.2, 0.5, 2, 5$, are reported in Tables 2–4. All results are based on 10,000 simulations.

Overall, the jackknife appears to be more accurate in terms of the size, especially when $m$ is relatively small (Case I). On the other hand, the simulated powers for POQUIM are higher at all alternatives, especially when $m$ is relatively small (Case I). However, it would be misleading to conclude that the POQUIM has higher power than the jackknife, because the power comparison is considered fair only if the two tests have similar sizes. In other words, for the case of $m = 50$, the higher power for POQUIM could be the result of the test overrejecting. Finally, note that the jackknife with the logarithmic transformation is specifically designed for this kind of model where the observations are divided into independent groups, while the POQUIM is for a much richer class of mixed linear models where the observations may or may not be divided into independent groups, as we will see in the next simulated example.

TABLE 1
*POQUIM versus jackknife—size*

| Nominal level | Method | Simulated size | | | | | |
|---|---|---|---|---|---|---|---|
| | | **I-i** | **I-ii** | **I-iii** | **II-i** | **II-ii** | **II-iii** |
| 0.01 | POQUIM | 0.022 | 0.026 | 0.028 | 0.011 | 0.013 | 0.015 |
| | Jackknife | 0.010 | 0.014 | 0.020 | 0.009 | 0.011 | 0.013 |
| 0.05 | POQUIM | 0.070 | 0.078 | 0.091 | 0.054 | 0.057 | 0.063 |
| | Jackknife | 0.052 | 0.053 | 0.068 | 0.053 | 0.053 | 0.060 |
| 0.10 | POQUIM | 0.123 | 0.132 | 0.151 | 0.106 | 0.108 | 0.114 |
| | Jackknife | 0.099 | 0.103 | 0.122 | 0.104 | 0.103 | 0.109 |

TABLE 2
*POQUIM versus jackknife—power (nominal level 0.01)*

| | | Simulated power | | | | | |
|---|---|---|---|---|---|---|---|
| Alternative | Method | I-i | I-ii | I-iii | II-i | II-ii | II-iii |
| $\gamma_1 = 0.2$ | POQUIM | 0.506 | 0.616 | 0.468 | 1.000 | 1.000 | 1.000 |
| | Jackknife | 0.487 | 0.463 | 0.454 | 1.000 | 1.000 | 1.000 |
| $\gamma_1 = 0.5$ | POQUIM | 0.112 | 0.164 | 0.122 | 0.914 | 0.891 | 0.793 |
| | Jackknife | 0.108 | 0.121 | 0.137 | 0.921 | 0.866 | 0.787 |
| $\gamma_1 = 2.0$ | POQUIM | 0.354 | 0.256 | 0.221 | 0.995 | 0.971 | 0.913 |
| | Jackknife | 0.196 | 0.118 | 0.072 | 0.993 | 0.968 | 0.887 |
| $\gamma_1 = 5.0$ | POQUIM | 0.991 | 0.954 | 0.900 | 1.000 | 1.000 | 1.000 |
| | Jackknife | 0.954 | 0.876 | 0.715 | 1.000 | 1.000 | 1.000 |

TABLE 3
*POQUIM versus jackknife—power (nominal level 0.05)*

| | | Simulated power | | | | | |
|---|---|---|---|---|---|---|---|
| Alternative | Method | I-i | I-ii | I-iii | II-i | II-ii | II-iii |
| $\gamma_1 = 0.2$ | POQUIM | 0.747 | 0.807 | 0.745 | 1.000 | 1.000 | 1.000 |
| | Jackknife | 0.728 | 0.709 | 0.668 | 1.000 | 1.000 | 1.000 |
| $\gamma_1 = 0.5$ | POQUIM | 0.283 | 0.336 | 0.286 | 0.980 | 0.966 | 0.917 |
| | Jackknife | 0.277 | 0.271 | 0.275 | 0.981 | 0.958 | 0.912 |
| $\gamma_1 = 2.0$ | POQUIM | 0.532 | 0.424 | 0.369 | 0.999 | 0.993 | 0.973 |
| | Jackknife | 0.411 | 0.317 | 0.223 | 0.999 | 0.993 | 0.970 |
| $\gamma_1 = 5.0$ | POQUIM | 0.997 | 0.984 | 0.956 | 1.000 | 1.000 | 1.000 |
| | Jackknife | 0.991 | 0.971 | 0.903 | 1.000 | 1.000 | 1.000 |

6.2. *A balanced two-way random effects model.* We now consider the balanced two-way random effects model of Example 1, also discussed in Section 3.1. Consider testing the hypothesis $H_0$: $\sigma_1^2 = \sigma_2^2$, or, equivalently, $H_0$: $\gamma_1 = \gamma_2$, which means that the two random effect factors contribute equally to the total variation. It is easy to show that the test statistic (27) reduces to

$$\hat{\chi}^2 = \frac{(\hat{\gamma}_1 - \hat{\gamma}_2)^2}{\hat{\Sigma}_{R,11} - 2\hat{\Sigma}_{R,12} + \hat{\Sigma}_{R,22}}, \tag{32}$$

where $\hat{\gamma}_1$ and $\hat{\gamma}_2$ are the REML estimators of $\gamma_1$ and $\gamma_2$, and $\hat{\Sigma}_{R,jk}$ is the $j,k$ element of the POQUIM estimator $\hat{\Sigma}_R$ of the ACM of $\hat{\theta} = (\hat{\lambda}, \hat{\gamma}_1, \hat{\gamma}_1)'$, the REML estimator. Note that in this case there are no (fully) specified values of the parameters under the null hypothesis, although the latter may still be used in some way (but the difference is expected to be small in large samples; see the remark below Theorem 4). On the other hand, it is interesting to see how the test performs when the straight POQUIM estimator is used in the denominator of (32), and that is what we do in this simulation. Once again, we study the performance of the test under both moderate and large sample sizes, as well as departures from normality. The following sample size configurations are considered: Case I, $m = 40$, $n = 40$; Case II, $m = 200$, $n = 200$. Furthermore, the following combinations of distributions for the random effects and errors are considered: Case i, $v$, $w \sim$ Normal; Case ii, $v$, $w \sim$ DE; Case iii, $v \sim$ DE, $w \sim$ CE; and Case iv, $v$, $w \sim$ CE. In all cases, $e \sim$ Normal. Note that the jackknife method discussed in the previous subsection does not apply to this case, because the observations cannot be divided into i.i.d. groups (or even independent groups). The true values of parameters are $\mu = \sigma_0^2 = \sigma_1^2 = 1.0$.

TABLE 4
*POQUIM versus jackknife—power (nominal level 0.10)*

| Alternative | Method | Simulated power | | | | | |
|---|---|---|---|---|---|---|---|
| | | **I-i** | **I-ii** | **I-iii** | **II-i** | **II-ii** | **II-iii** |
| $\gamma_1 = 0.2$ | POQUIM | 0.844 | 0.875 | 0.807 | 1.000 | 1.000 | 1.000 |
| | Jackknife | 0.829 | 0.810 | 0.776 | 1.000 | 1.000 | 1.000 |
| $\gamma_1 = 0.5$ | POQUIM | 0.405 | 0.442 | 0.396 | 0.991 | 0.983 | 0.954 |
| | Jackknife | 0.398 | 0.382 | 0.372 | 0.991 | 0.979 | 0.950 |
| $\gamma_1 = 2.0$ | POQUIM | 0.633 | 0.564 | 0.462 | 1.000 | 0.997 | 0.987 |
| | Jackknife | 0.540 | 0.453 | 0.350 | 1.000 | 0.997 | 0.986 |
| $\gamma_1 = 5.0$ | POQUIM | 0.999 | 0.992 | 0.975 | 1.000 | 1.000 | 1.000 |
| | Jackknife | 0.998 | 0.988 | 0.954 | 1.000 | 1.000 | 1.000 |

As in the previous subsection, we first consider the size of the test, so we take $\sigma_2^2 = 1.0$. The simulated sizes corresponding to the nominal levels 0.01, 0.05 and 0.10 are reported in Table 5. Next we look at the powers at the following alternatives: $\sigma_2^2 = 0.2$, 0.5, 2, 5, which correspond to $\gamma_2/\gamma_1 = 0.2$, 0.5, 2, 5, respectively. The simulated powers are reported in Tables 6–8. Again, all results are based on 10,000 simulations.

The numbers seem to follow the same pattern. As the sample size increases, the simulated sizes get closer to the nominal levels and the simulated powers increase significantly. There does not seem to be a difference, in terms of the size, across different distributions. However, the simulated powers appear significantly higher when all the distributions are normal as compared to other cases where the distributions of the random effects are nonnormal. Also, the powers are relatively low when the alternatives are close to the null ($\gamma_2/\gamma_1 = 0.5$ or 2.0), but much improved when the alternatives are further away ($\gamma_2/\gamma_1 = 0.2$ and 5.0). Overall, the simulation results are consistent with the theoretical findings of Theorem 4.

TABLE 5
*Simulated size*

| Nominal level | I-i | I-ii | I-iii | I-iv | II-i | II-ii | II-iii | II-iv |
|---|---|---|---|---|---|---|---|---|
| 0.01 | 0.014 | 0.011 | 0.014 | 0.011 | 0.011 | 0.008 | 0.011 | 0.008 |
| 0.05 | 0.071 | 0.061 | 0.070 | 0.066 | 0.053 | 0.051 | 0.055 | 0.048 |
| 0.10 | 0.135 | 0.126 | 0.139 | 0.136 | 0.108 | 0.108 | 0.109 | 0.102 |

TABLE 6
*Simulated power (nominal level 0.01)*

| Alternative | I-i | I-ii | I-iii | I-iv | II-i | II-ii | II-iii | II-iv |
|---|---|---|---|---|---|---|---|---|
| $\gamma_2/\gamma_1 = 0.2$ | 0.955 | 0.568 | 0.551 | 0.398 | 1.000 | 1.000 | 0.999 | 0.986 |
| $\gamma_2/\gamma_1 = 0.5$ | 0.313 | 0.100 | 0.118 | 0.073 | 0.988 | 0.684 | 0.619 | 0.439 |
| $\gamma_2/\gamma_1 = 2.0$ | 0.324 | 0.100 | 0.070 | 0.088 | 0.988 | 0.685 | 0.459 | 0.443 |
| $\gamma_2/\gamma_1 = 5.0$ | 0.969 | 0.649 | 0.491 | 0.497 | 1.000 | 0.999 | 0.989 | 0.992 |

TABLE 7
*Simulated power (nominal level 0.05)*

| Alternative | I-i | I-ii | I-iii | I-iv | II-i | II-ii | II-iii | II-iv |
|---|---|---|---|---|---|---|---|---|
| $\gamma_2/\gamma_1 = 0.2$ | 0.994 | 0.864 | 0.839 | 0.713 | 1.000 | 1.000 | 1.000 | 0.999 |
| $\gamma_2/\gamma_1 = 0.5$ | 0.579 | 0.308 | 0.321 | 0.232 | 0.998 | 0.874 | 0.819 | 0.713 |
| $\gamma_2/\gamma_1 = 2.0$ | 0.595 | 0.305 | 0.227 | 0.256 | 0.998 | 0.879 | 0.764 | 0.717 |
| $\gamma_2/\gamma_1 = 5.0$ | 0.997 | 0.901 | 0.799 | 0.779 | 1.000 | 1.000 | 1.000 | 0.999 |

TABLE 8
*Simulated power (nominal level 0.10)*

| Alternative | I-i | I-ii | I-iii | I-iv | II-i | II-ii | II-iii | II-iv |
|---|---|---|---|---|---|---|---|---|
| $\gamma_2/\gamma_1 = 0.2$ | 0.999 | 0.946 | 0.923 | 0.846 | 1.000 | 1.000 | 1.000 | 1.000 |
| $\gamma_2/\gamma_1 = 0.5$ | 0.702 | 0.456 | 0.451 | 0.364 | 0.999 | 0.931 | 0.887 | 0.818 |
| $\gamma_2/\gamma_1 = 2.0$ | 0.719 | 0.448 | 0.359 | 0.382 | 0.999 | 0.936 | 0.864 | 0.818 |
| $\gamma_2/\gamma_1 = 5.0$ | 0.999 | 0.955 | 0.904 | 0.879 | 1.000 | 1.000 | 1.000 | 1.000 |

**7. Discussion and remarks.** A classic parametric statistical model assumes that the distribution of the data is fully determined by a vector of parameters. Under such a model, a maximum likelihood estimator of the vector of parameters is self-contained in the sense that the (asymptotic) covariance matrix of the estimator does not involve any additional unknown parameter. In many cases, however, a model is not fully determined by a set of parameters. For example, under nonnormality, the distribution of the data is not determined by the mean and the variance. Obviously, in such cases maximum likelihood does not apply, but a quasi-likelihood method may be used to estimate the parameters of direct interest (e.g., [11]). The problem is that the estimator may no longer be self-contained. The POQUIM method provides a way to estimate the (asymptotic) covariance matrix of a maximum quasi-likelihood estimator and, hence, self-contains the latter.

The general procedure of POQUIM is the following: Let $l(\theta)$ be the quasi-log-likelihood. Then the ACM of $\hat{\theta}$, the maximum quasi-likelihood estimator, is $\Sigma = \mathcal{I}_2^{-1}\mathcal{I}_1\mathcal{I}_2^{-1}$, where $\mathcal{I}_1 = \text{Var}(\partial l/\partial \theta)$ and $\mathcal{I}_2 = \text{E}(\partial^2 l/\partial \theta\, \partial \theta')$. Usually, $\mathcal{I}_2$ either does not involve additional parameters or, if it does, at least it can be estimated by an observed form (e.g., by $\partial^2 l/\partial \theta\, \partial \theta'$ with $\theta$ replaced by $\hat{\theta}$). However, when the data are correlated, the matrix $\mathcal{I}_1$ cannot be estimated by an observed form. The idea of POQUIM is to express $\mathcal{I}_1$ as $\text{E}(S_1) + S_2$ such that $\text{E}(S_1)$ involves parameters other than $\theta$ but that can be estimated by an observed form (i.e., by $S_1$ with $\theta$ replaced by $\hat{\theta}$), and $S_2$ does not involve any additional parameters and therefore can be estimated by an estimated form (i.e., by $S_2$ with $\theta$ replaced by $\hat{\theta}$). In this paper this general method is applied to mixed linear models.

In this paper we assume that the random effects and errors in a non-Gaussian mixed linear model are independent. This means that (i) the vectors $\alpha_1, \ldots, \alpha_s, \varepsilon$ are independent and (ii) the components of $\alpha_j$ $(1 \le j \le s)$ and $\varepsilon$ are independent. A mixed linear model that satisfies (1) as well as (i) and (ii) above is also known as an analysis of variance mixed (ANOVA) model (e.g., [5], where normality is assumed). Mixed ANOVA models are very popular in practice. On the other hand, there are also mixed linear models used in practice that involve dependent random effects or errors, such as the so-called longitudinal model (e.g., [5], [18]). For example, in Section 3.4 it may be reasonable to assume that the random intercept, $a_i$, and slope, $b_i$, which correspond to the same individual, are correlated. Note that, in cases of correlated random effects, the variance components are defined as the parameters involved in $V$, the covariance matrix of $y$, that involve correlations in addition to the variances. Furthermore, there are more additional parameters involved in the ACM of, say, the REML estimator. However, a possible POQUIM decomposition may still be obtained. For example, by (9), one can write

$$\text{cov}\left(\frac{\partial l_\text{R}}{\partial \theta_j}, \frac{\partial l_\text{R}}{\partial \theta_k}\right) = \sum_{(i_1,\ldots,i_4)\in\mathcal{S}_1} B_{j,i_1,i_2}B_{k,i_3,i_4}\,\text{cov}(u_{i_1}u_{i_2}, u_{i_3}u_{i_4})$$
$$+ \sum_{(i_1,\ldots,i_4)\in\mathcal{S}_2} B_{j,i_1,i_2}B_{k,i_3,i_4}\,\text{cov}(u_{i_1}u_{i_2}, u_{i_3}u_{i_4})$$
$$= I_1 + I_2,$$

where $\mathcal{S}_2$ are those indices such that $\text{cov}(u_{i_1}u_{i_2}, u_{i_3}u_{i_4})$ involve only the variance components and $\mathcal{S}_1$ are those indexes such that $\text{cov}(u_{i_1}u_{i_2}, u_{i_3}u_{i_4})$ involve additional parameters. If $I_1$ can be further expressed as an expected value plus a term that depends only on the variance components, one has a potential POQUIM decomposition.

The robust dispersion test derived in Section 5 is of $\chi^2$ type. As mentioned in Section 6.1, $\chi^2$ tests are *omnibus* (e.g., [21], Section 1.1). However, a *directional* test could be obtained, using asymptotic normality of, say, the REML estimator $\hat{\theta}$ and the POQUIM estimator of its ACM. The only exception is when $\theta$ lies on the boundary of the parameter space, because, obviously, in this case $\hat{\theta}$ cannot be asymptotically normal if it is required to stay in the parameter space. However, for *testing purposes* one may relax the latter restriction, for example, by defining $\hat{\theta}$ as the solution to the REML equation, which may or may not be in the parameter space (see the remark above Theorem 2). With such a definition, asymptotic normality of $\hat{\theta}$ may still hold, even if $\theta$ is on the boundary of the parameter space. See, for example, [22]. Such a result may be used for testing, for example, that some of the variance components are zero.

The conditions of Theorem 2 [more specifically, condition (iii)] imply the existence and consistency of the REML estimator. See the remarks below Theorem 2 and also those above Theorem 2 regarding the definition of the REML estimator. Of course, this is a large sample result, which does not guarantee the existence of the REML estimator in a finite sample situation, even under the normality assumption. A similar problem exists for the ML estimator as well. See, for example, [6] and [4].

## 8. Proofs and other technical details.

### 8.1. *Partial derivatives of the quasi-restricted log-likelihood.* Differentiating (5) with respect to $\theta$ and using the fact $PVP = P$, we have $\partial l_{\mathrm{R}}/\partial\lambda = (y'Py - N + p)/2\lambda$ and $\partial l_{\mathrm{R}}/\partial\gamma_j = (\lambda/2)\{y'PZ_jZ_j'Py - \mathrm{tr}(PZ_jZ_j')\}$, $1 \le j \le s$. Note that since $PX = 0$, the vector $y$ in the above expressions can be replaced by $u = y - X\beta$. Furthermore, we have $\mathrm{E}(\partial^2 l_{\mathrm{R}}/\partial\lambda^2) = -(N - p)/2\lambda^2$, $\mathrm{E}(\partial^2 l_{\mathrm{R}}/\partial\lambda\,\partial\gamma_j) = -(1/2)\mathrm{tr}(PZ_jZ_j')$, $1 \le j \le s$, and $\mathrm{E}(\partial^2 l_{\mathrm{R}}/\partial\gamma_j\,\partial\gamma_k) = -(\lambda^2/2)\mathrm{tr}(PZ_jZ_j'PZ_kZ_k')$, $1 \le j, k \le s$.

### 8.2. *Proof of Lemma 1.* First note $u_i = \sum_{t=0}^s u_{it}$ with $u_{it} = \sum_{l=1}^{m_t} z_{itl}\alpha_{tl}$; hence $\mathrm{E}(u_{i_1}u_{i_2}) = \sum_{t=0}^s \mathrm{E}(u_{i_1t}u_{i_2t}) = \sum_{t=0}^s \sigma_t^2 z_{i_1t}' z_{i_2t} = \lambda\Gamma(i_1, i_2)$. Next, we have $\mathrm{E}(u_{i_1}\cdots u_{i_4}) = \sum_{t_1,\dots,t_4}\mathrm{E}(u_{i_1t_1}\cdots u_{i_4t_4})$ and $\mathrm{E}(u_{i_1t_1}\cdots u_{i_4t_4}) = 0$ unless (1) $t_1 = \cdots = t_4$ or (2) the $t$'s are in two pairs. Furthermore, under (1) we have $\mathrm{E}(u_{i_1t_1}\cdots u_{i_4t_4}) = \sum_{l_1,\dots,l_4} z_{i_1tl_1}\cdots z_{i_4tl_4}\mathrm{E}(\alpha_{tl_1}\cdots\alpha_{tl_4})$, where $t_1 = \cdots = t_4 = t$. Again, $\mathrm{E}(\alpha_{tl_1}\cdots\alpha_{tl_4}) = 0$ unless (1–1) $l_1 = \cdots = l_4$, (1–2) $l_1 = l_2 \ne l_3 = l_4$, (1–3) $l_1 = l_3 \ne l_2 = l_4$ or (1–4) $l_1 = l_4 \ne l_2 = l_3$. It is easy to show that $\sum_{(1-1)}\cdots = \mathrm{E}(\alpha_{t1}^4)z_{i_1t}\cdots z_{i_4t}$, $\sum_{(1-2)}\cdots = \sigma_t^4\{(z_{i_1t}\cdot z_{i_2t})(z_{i_3t}\cdot$

$z_{i_4t}) - z_{i_1t} \cdots z_{i_4t}\}$, $\sum_{(1-3)} \cdots = \sigma_t^4\{(z_{i_1t} \cdot z_{i_3t})(z_{i_2t} \cdot z_{i_4t}) - z_{i_1t} \cdots z_{i_4t}\}$ and $\sum_{(1-4)} \cdots = \sigma_t^4\{(z_{i_1t} \cdot z_{i_4t})(z_{i_2t} \cdot z_{i_3t}) - z_{i_1t} \cdots z_{i_4t}\}$. It follows that

$$\sum_{(1)} \mathrm{E}(u_{i_1t_1} \cdots u_{i_4t_4})$$

$$= \sum_t \kappa_t z_{i_1t} \cdots z_{i_4t}$$

$$+ \lambda^2 \Bigg\{ \sum_t \gamma_t^2 (z_{i_1t} \cdot z_{i_2t})(z_{i_3t} \cdot z_{i_4t})$$

$$+ \sum_t \gamma_t^2 (z_{i_1t} \cdot z_{i_3t})(z_{i_2t} \cdot z_{i_4t}) + \sum_t \gamma_t^2 (z_{i_1t} \cdot z_{i_4t})(z_{i_2t} \cdot z_{i_3t}) \Bigg\}.$$

Similarly, (2) has three cases: (2–1) $t_1 = t_2 \neq t_3 = t_4$, (2–2) $t_1 = t_3 \neq t_2 = t_3$ and (2–3) $t_1 = t_4 \neq t_2 = t_3$. Furthermore, we have

$$\sum_{(2-1)} \mathrm{E}(u_{i_1t_1} \cdots u_{i_4t_4}) = \lambda^2 \Bigg\{ \Gamma(i_1, i_2)\Gamma(i_3, i_4) - \sum_t \gamma_t^2 (z_{i_1t} \cdot z_{i_2t})(z_{i_3t} \cdot z_{i_4t}) \Bigg\},$$

$$\sum_{(2-2)} \mathrm{E}(u_{i_1t_1} \cdots u_{i_4t_4}) = \lambda^2 \Bigg\{ \Gamma(i_1, i_3)\Gamma(i_2, i_4) - \sum_t \gamma_t^2 (z_{i_1t} \cdot z_{i_3t})(z_{i_2t} \cdot z_{i_4t}) \Bigg\},$$

$$\sum_{(2-3)} \mathrm{E}(u_{i_1t_1} \cdots u_{i_4t_4}) = \lambda^2 \Bigg\{ \Gamma(i_1, i_4)\Gamma(i_2, i_3) - \sum_t \gamma_t^2 (z_{i_1t} \cdot z_{i_4t})(z_{i_2t} \cdot z_{i_3t}) \Bigg\}.$$

Therefore, in conclusion, we have $\mathrm{cov}(u_{i_1}u_{i_2}, u_{i_3}u_{i_4}) = \mathrm{E}(u_{i_1} \cdots u_{i_4}) - \lambda^2 \Gamma(i_1, i_2)\Gamma(i_3, i_4)$ equal to the right-hand side of (6).

Equation (7) is easily derived from (6), observing the first equation of (9).

8.3. *Proof of Theorem* 2. In the following discussion $A_N, \ldots$ represent sequences of matrices (vectors, numbers), but for notational simplicity we suppress the subscript $N$. Note that for a matrix $B$, $B > 0$ means that $B$ is positive definite. We first state and prove a lemma.

LEMMA 3. *Let $A$, $G$ be sequences of positive definite matrices such that $G^{-1}AG^{-1} \to B > 0$. Let $\tilde{A}$ be another sequence of matrices. Then $A^{-1/2}\tilde{A} \times A^{-1/2} \to I$, the identity matrix, if and only if $G^{-1}(\tilde{A} - A)G^{-1} \to 0$.*

PROOF. Suppose that $A^{-1/2}\tilde{A}A^{-1/2} \to I$. Then we have

$$G^{-1}(\tilde{A} - A)G^{-1} = G^{-1}A^{1/2}(A^{-1/2}\tilde{A}A^{-1/2} - I)A^{1/2}G^{-1}.$$

Therefore,

$$\|G^{-1}(\tilde{A} - A)G^{-1}\| \leq \|G^{-1}A^{1/2}\|^2\|A^{-1/2}\tilde{A}A^{-1/2} - I\|$$
$$= \lambda_{\max}(G^{-1}AG^{-1})\|A^{-1/2}\tilde{A}A^{-1/2} - I\| \longrightarrow 0.$$

Now suppose that $G^{-1}(\tilde{A} - A)G^{-1} \to 0$. Let $\Delta = G^{-1}AG^{-1} - B$. Then we have

$$A^{-1/2}\tilde{A}A^{-1/2} - I = A^{-1/2}GBGA^{-1/2} - I$$
$$+ A^{-1/2}GB^{1/2}(B^{-1/2}G^{-1}\tilde{A}G^{-1}B^{-1/2} - I)B^{1/2}GA^{-1/2}$$
$$= D_1 + D_2.$$

We have $D_1 = A^{-1/2}(GBG - A)A^{-1/2} = -A^{-1/2}G\Delta GA^{-1/2}$; thus $\|D_1\| \leq \|A^{-1/2}G\|^2\|\Delta\| = \lambda_{\max}(GA^{-1}G)\|\Delta\| \to 0$, because $GA^{-1}G \to B^{-1}$ and $\Delta \to 0$. On the other hand, we have

$$B^{-1/2}G^{-1}\tilde{A}G^{-1}B^{-1/2} - I = B^{-1/2}G^{-1}AG^{-1}B^{-1/2} - I$$
$$+ B^{-1/2}G^{-1}(\tilde{A} - A)G^{-1}B^{-1/2} \longrightarrow 0.$$

Therefore,

$$\|D_2\| \leq \|A^{-1/2}GB^{1/2}\|^2\|B^{-1/2}G^{-1}\tilde{A}G^{-1}B^{-1/2} - I\|$$
$$= \lambda_{\max}(B^{1/2}GA^{-1}GB^{1/2})\|B^{-1/2}G^{-1}\tilde{A}G^{-1}B^{-1/2} - I\| \longrightarrow 0. \quad \square$$

Recall that a sequence of matrices $M$ is bounded from above if $\|M\|$ is bounded; the sequence is bounded from below if $\|M^{-1}\|$ is bounded.

COROLLARY 1. *Let $A$, $G$ be sequences of positive definite matrices such that $G^{-1}AG^{-1}$ is bounded from above as well as from below. Let $\hat{A}$ be a sequence of random matrices. Then $A^{-1/2}\hat{A}A^{-1/2} \to I$ in probability if and only if $G^{-1}(\hat{A} - A)G^{-1} \to 0$ in probability.*

PROOF. A sequence of random matrices converges in probability if and only if for any subsequence, there is a further subsequence that converges almost surely (to the same limit).

First assume that $A^{-1/2}\hat{A}A^{-1/2} \to I$ in probability. For any subsequence of $G^{-1}(\hat{A} - A)G^{-1}$, since the corresponding subsequence of $G^{-1}AG^{-1}$ is bounded, there is a further subsequence such that $G^{-1}AG^{-1} \to B$ for some $B > 0$. The latter property is implied by the boundedness from below of the subsequence. Consider the corresponding further subsequence of $A^{-1/2}\hat{A}A^{-1/2}$. Since it converges in probability, there is a further subsequence such that $A^{-1/2}\hat{A}A^{-1/2} \to I$ almost surely. It follows, by Lemma 1, that the corresponding further subsequence $G^{-1}(\hat{A} - A)G^{-1} \to 0$ almost surely.

Next assume that $G^{-1}(\hat{A} - A)G^{-1} \to 0$ in probability. By similar arguments we have for any subsequence of $A^{-1/2}\hat{A}A^{-1/2}$ that there is a further subsequence that $\to I$ almost surely. $\square$

The next lemma states that the Gaussian information matrix $\mathcal{G}$ and the QUIM $\mathcal{I}_1$ are asymptotically of the same order.

LEMMA 4. *Under condition* (i) *of Theorem* 2, *there are positive constants $a$ and $b$ such that $a\mathcal{G} \leq \mathcal{I}_1 \leq b\mathcal{G}$.*

PROOF. First, it is easy to show that

$$\text{tr}(B_j V B_k V) = \sum_{t_1, t_2 = 0}^{s} \sigma_{t_1}^2 \sigma_{t_2}^2 \sum_{l_1 = 1}^{m_{t_1}} \sum_{l_2 = 1}^{m_{t_2}} (z'_{t_1 l_1} B_j z_{t_2 l_2})(z'_{t_1 l_1} B_k z_{t_2 l_2}).$$

Condition (i) of Theorem 2 implies that there is $0 < \delta < 1$ such that $\kappa_t = \text{var}(\alpha_{t1}^2) - 2\sigma_t^4 \geq 2(\delta - 1)\sigma_t^4$, $0 \leq t \leq s$. Thus, it can be shown by (7) that for any $x = (x_j)_{0 \leq j \leq s}$,

$$x'\mathcal{I}_1 x = x'\mathcal{G}x + \sum_{t=0}^{s} \kappa_t \sum_{l=1}^{m_t} \left\{ \sum_{j=0}^{s} x_j (z'_{tl} B_j z_{tl}) \right\}^2$$

$$\geq \delta x'\mathcal{G}x + 2(1-\delta)\left[ \sum_{t_1, t_2 = 0}^{s} \sigma_{t_1}^2 \sigma_{t_2}^2 \sum_{l_1 = 1}^{m_{t_1}} \sum_{l_2 = 1}^{m_{t_2}} \left\{ \sum_{j=0}^{s} x_j (z'_{t_1 l_1} B_j z_{t_2 l_2}) \right\}^2 \right.$$

$$\left. - \sum_{t=0}^{s} \sigma_t^4 \sum_{l=1}^{m_t} \left\{ \sum_{j=0}^{s} x_j (z'_{tl} B_j z_{tl}) \right\}^2 \right]$$

$$\geq \delta x'\mathcal{G}x.$$

On the other hand, condition (i) of Theorem 2 implies that there is $M > 0$ such that $\kappa_4 \leq 2M\sigma_t^4$, $0 \leq t \leq s$. Thus, we have, similarly,

$$x'\mathcal{I}_1 x \leq x'\mathcal{G}x + 2M \sum_{t=0}^{s} \sigma_t^4 \sum_{l=1}^{m_t} \left\{ \sum_{j=0}^{s} x_j (z'_{tl} B_j z_{tl}) \right\}^2$$

$$\leq (1+M)x'\mathcal{G}x. \qquad \square$$

COROLLARY 2. *Under condition* (i) *of Theorem* 2, *$G^{-1}\mathcal{G}G^{-1}$ is bounded from above as well as from below if and only if $G^{-1}\mathcal{I}_1 G^{-1}$ is bounded from above as well as from below.*

Throughout the rest of the proof, $c$ represents a positive constant whose value may be different in each occurrence.

First prove that $\hat{\mathcal{I}}_1$ is consistent. According to condition (iii) and Corollary 2, the sequence $G^{-1}\mathcal{I}_1 G^{-1}$ is bounded from above as well as from below. Then, according to Corollary 1, it suffices to show that $G^{-1}(\hat{\mathcal{I}}_1 - \mathcal{I}_1)G^{-1} \to 0$ in probability.

First consider the observed part. Let $\mathcal{D} = \{|\hat{\theta} - \theta| \le \delta, |\hat{\beta} - \beta| \le \delta\}$ $(\delta > 0)$ and $\tilde{\mathcal{I}}_{1,1,jk} = \sum_{l=1}^{L} c_{j,k,l} \sum_{f(i_1,\dots,i_4)=f_l} u_{i_1} \cdots u_{i_4}$. It is easy to show that, on $\mathcal{D}$, $|\hat{u}_{i_1} \cdots \hat{u}_{i_4} - u_{i_1} \cdots u_{i_4}| \le c\delta \sum_{r=1}^{4}(y_{i_r}^4 + |x_{i_r}|^4)$ and $|\hat{u}_{i_1} \cdots \hat{u}_{i_4}| \le c \sum_{r=1}^{4}(y_{i_r}^4 + |x_{i_r}|^4)$. It follows that, on $\mathcal{D}$,

$$\left| \hat{c}_{j,k,l} \sum_{f(i_1,\dots,i_4)=f_l} \hat{u}_{i_1} \cdots \hat{u}_{i_4} - c_{j,k,l} \sum_{f(i_1,\dots,i_4)=f_l} u_{i_1} \cdots u_{i_4} \right|$$

$$\le c\{d_{j,k,l}(\delta) + |c_{j,k,l}|\delta\} \sum_{f(i_1,\dots,i_4)=f_l} \sum_{r=1}^{4}(y_{i_r}^4 + |x_{i_r}|^4).$$

Conditions (i) and (ii) imply that $\mathrm{E}(y_t^4) \le c$. Therefore, we have

$$\begin{aligned}
&\mathrm{E}\{(g_j g_k)^{-1}|\hat{\mathcal{I}}_{1,1,jk} - \tilde{\mathcal{I}}_{1,1,jk}|\mathbf{1}_{\mathcal{D}}\} \\
(33)\quad &\le c(g_j g_k)^{-1} \sum_{l=1}^{L}\{d_{j,k,l}(\delta) + |c_{j,k,l}|\delta\} \sum_{f(i_1,\dots,i_4)=f_l} \sum_{r=1}^{4}\{\mathrm{E}(y_{i_r}^4) + |x_{i_r}|^4\} \\
&\le c\left\{ (g_j g_k)^{-1} \sum_{l=1}^{L} h_l d_{j,k,l}(\delta) + \delta(g_j g_k)^{-1} \sum_{l=1}^{L} h_l |c_{j,k,l}| \right\}.
\end{aligned}$$

On the other hand, given $M > 0$, we have $\alpha_{tl} = \alpha_{tl1} + \alpha_{tl2}$, where $\alpha_{tl1} = \alpha_{tl}\mathbf{1}_{(|\alpha_{tl}| \le M)} - \mathrm{E}\{\alpha_{tl}\mathbf{1}_{(|\alpha_{tl}| \le M)}\}$. Thus,

$$u_i = \sum_{t=0}^{s} z_{it}' \alpha_t = u_{i1} + u_{i2},$$

where $u_{ir} = \sum_{t=0}^{s} \sum_{l=1}^{m_t} z_{itl}\alpha_{tlr}$, $r = 1, 2$. Conditions (i) and (ii) imply that $\mathrm{E}(u_{i1}^4) \le c \sum_{t=0}^{s} \|z_{it}\|_1^4 \le c$ and $\mathrm{E}(u_{i2}^4) \le cb^4(M)$, where $b(M) = \sum_{t=0}^{s} \mathrm{E}(\alpha_{t1}^4 \mathbf{1}_{(|\alpha_{t1}| > M)}) \to 0$ as $M \to \infty$. Write

$$u_{i_1} \cdots u_{i_4} - \mathrm{E}(u_{i_1} \cdots u_{i_4}) = u_{i_1 1} \cdots u_{i_4 1} - \mathrm{E}(u_{i_1 1} \cdots u_{i_4 1}) + (\cdots) - \mathrm{E}(\cdots),$$

where $(\cdots)$ is a sum of products with each product involving at least one $u_{i_r 2}$ $(r = 1, \dots, 4)$. It follows by Hölder's inequality that $|\mathrm{E}(\cdots)| \le \mathrm{E}|\cdots| \le cb(M)$. Therefore, we can write

$$\tilde{\mathcal{I}}_{1,1,jk} - \mathcal{I}_{1,1,jk} = \sum_{l=1}^{L} c_{j,k,l} \sum_{f(i_1,\dots,i_4)=f_l} \{u_{i_1 1} \cdots u_{i_4 1} - \mathrm{E}(u_{i_1 1} \cdots u_{i_4 1})\}$$

$$+ \sum_{l=1}^{L} c_{j,k,l} \sum_{f(i_1,\ldots,i_4)=f_l} \{(\cdots) - \mathrm{E}(\cdots)\}$$

$$= S_1 + S_2$$

with $\mathrm{E}(|S_2|/g_j g_k) \le c b(M)(g_j g_k)^{-1} \sum_{l=1}^{L} h_l |c_{j,k,l}|$. Furthermore, we have

$$\mathrm{E}(S_1^2) = \sum_{l_1, l_2 = 1}^{L} c_{j,k,l_1} c_{j,k,l_2}$$

$$\times \sum_{f(i_1,\ldots,i_4)=f_{l_1}, f(i_5,\ldots,i_8)=f_{l_2}} \mathrm{cov}(u_{i_1 1} \cdots u_{i_4 1}, u_{i_5 1} \cdots u_{i_8 1}).$$

The nonzero components of $d_{i_1} + \cdots + d_{i_4}$ ($d_i$ is defined in the second paragraph of Section 2.2) correspond to the indexes of the random effects and errors involved in $u_{i_1 1} \cdots u_{i_4 1}$. Thus, if $(d_{i_1} + \cdots + d_{i_4}) \cdot (d_{i_5} + \cdots + d_{i_8}) = 0$, $u_{i_1 1} \cdots u_{i_4 1}$ and $u_{i_5 1} \cdots u_{i_8 1}$ involve different random effects and errors, hence $\mathrm{cov}(u_{i_1 1} \cdots u_{i_4 1}, u_{i_5 1} \cdots u_{i_8 1}) = 0$; otherwise, the covariance is bounded in absolute value by $cM^8$. It follows that

$$\mathrm{E}(S_1/g_j g_k)^2 \le cM^8 (g_j g_k)^{-2} \sum_{l_1, l_2 = 1}^{L} h_{l_1, l_2} |c_{j,k,l_1} c_{j,k,l_2}|.$$

In conclusion, we have

$$\mathrm{E}\{(g_j g_k)^{-1} |\tilde{\mathcal{I}}_{1,1,jk} - \mathcal{I}_{1,1,jk}|\}$$

$$(34) \qquad \le c b(M)(g_j g_k)^{-1} \sum_{l=1}^{L} h_l |c_{j,k,l}|$$

$$+ cM^4 \sqrt{(g_j g_k)^{-2} \sum_{l_1, l_2 = 1}^{L} h_{l_1, l_2} |c_{j,k,l_1} c_{j,k,l_2}|}.$$

Now consider the estimated part. We have

$$\hat{\mathcal{I}}_{1,2,jk} - \mathcal{I}_{1,2,jk}$$

$$= 2\{\mathrm{tr}(\hat{B}_j \hat{V} \hat{B}_k \hat{V}) - \mathrm{tr}(B_j V B_k V)\}$$

$$- 3\hat{\lambda}^2 \sum_{l=1}^{L} (\hat{c}_{j,k,l} - c_{j,k,l}) \sum_{f(i_1,\ldots,i_4)=f_l} \hat{\Gamma}(i_1, i_3) \hat{\Gamma}(i_2, i_4)$$

$$- 3\hat{\lambda}^2 \sum_{l=1}^{L} c_{j,k,l} \sum_{f(i_1,\ldots,i_4)=f_l} \{\hat{\Gamma}(i_1, i_3) \hat{\Gamma}(i_2, i_4) - \Gamma(i_1, i_3) \Gamma(i_2, i_4)\}$$

$$-3(\hat{\lambda}^2 - \lambda^2) \sum_{l=1}^{L} c_{j,k,l} \sum_{f(i_1,\ldots,i_4)=f_l} \Gamma(i_1,i_3)\Gamma(i_2,i_4)$$

$$= T_1 - \sum_{r=1}^{3} T_r.$$

On $\mathcal{D}$ we have $|T_1| \leq g_{j,k}(\delta)$. Also, condition (ii) implies that, on $\mathcal{D}$, $|T_2| \leq c \sum_{l=1}^{L} h_l d_{j,k,l}(\delta)$ and $|T_r| \leq c\delta \sum_{l=1}^{L} h_l |c_{j,k,l}|$, $r = 3, 4$. Therefore, we have

$$(g_j g_k)^{-1}|\hat{\mathcal{I}}_{1,2,jk} - \mathcal{I}_{1,2,jk}|$$

$$(35) \qquad \leq (g_j g_k)^{-1} g_{j,k}(\delta) + c(g_j g_k)^{-1} \sum_{l=1}^{L} h_l d_{j,k,l}(\delta)$$

$$+ c\delta(g_j g_k)^{-1} \sum_{l=1}^{L} h_l |c_{j,k,l}| \qquad \text{on } \mathcal{D}.$$

For any $\eta > 0$ and $\rho > 0$, first choose $\delta > 0$ and $M > 0$ such that the right-hand side of (33) is less than $\eta\rho/15$, the right-hand side of (35) is less than $\eta/3$ and the first term on the right-hand side of (34) is less than $\eta\rho/15$, which one can do according to conditions (i), (ii), (iv) and (v). Then choose $N_0$ such that, when $N \geq N_0$, we have $P(\mathcal{D}^c) < \rho/5$, and the second term on the right-hand side of (34) is less than $\eta\rho/15$, which one can do by conditions (iii) and (iv). Note that condition (iii) implies consistency of $\hat{\theta}$ and $\hat{\beta}$ (see the remarks below Theorem 2). It follows that, when $N \geq N_0$, by (33) and Chebyshev's inequality,

$$P\left((g_j g_k)^{-1}|\hat{\mathcal{I}}_{1,1,jk} - \tilde{\mathcal{I}}_{1,1,jk}| > \frac{\eta}{3}\right)$$

$$\leq P\left((g_j g_k)^{-1}|\hat{\mathcal{I}}_{1,1,jk} - \tilde{\mathcal{I}}_{1,1,jk}|\mathbf{1}_{\mathcal{D}} > \frac{\eta}{3}\right) + P(\mathcal{D}^c)$$

$$\leq \frac{3}{\eta} E\{(g_j g_k)^{-1}|\hat{\mathcal{I}}_{1,1,jk} - \tilde{\mathcal{I}}_{1,1,jk}|\mathbf{1}_{\mathcal{D}}\} + \frac{\rho}{5} < \frac{2}{5}\rho.$$

Similarly, by (34), $P((g_j g_k)^{-1}|\tilde{\mathcal{I}}_{1,1,jk} - \mathcal{I}_{1,1,jk}| > \eta/3) < (2/5)\rho$, and by (35), $P((g_j g_k)^{-1}|\hat{\mathcal{I}}_{1,2,jk} - \mathcal{I}_{1,2,jk}| > \eta/3) \leq P(\mathcal{D}^c) < \rho/5$. Therefore, we conclude that when $N \geq N_0$, the probability is greater than $1 - \rho$ that $(g_j g_k)^{-1} \times |\hat{\mathcal{I}}_{1,jk} - \mathcal{I}_{1,jk}| \leq \eta$.

Now consider consistency of $\hat{\Sigma}_R$. First note that $\mathcal{I}_2 = -\mathcal{G}$, which is nonsingular by condition (iii). Since $G\Sigma_R G = (G^{-1}\mathcal{G}G^{-1})^{-1} G^{-1} \mathcal{I}_1 G^{-1} (G^{-1}\mathcal{G}G^{-1})^{-1}$, by Corollary 2, $G\Sigma_R G$ is bounded from above as well as from below. Then, by Corollary 1 [note that $G = (G^{-1})^{-1}$], it suffices to show that $G(\hat{\Sigma}_R - \Sigma_R)G \to 0$ in probability.

By condition (v) and consistency of $\hat{\theta}$, it is easy to show that $G^{-1}(\hat{\mathcal{G}} - \mathcal{G})G^{-1} \to 0$ in probability, where $\hat{\mathcal{G}}$ is $\mathcal{G}$ with $\theta$ replaced by $\hat{\theta}$. Note that $\hat{\mathcal{I}}_2 = -\hat{\mathcal{G}}$. Also, since $G^{-1}\hat{\mathcal{G}}G^{-1} = G^{-1}\mathcal{G}G^{-1} + G^{-1}(\hat{\mathcal{G}} - \mathcal{G})G^{-1}$, we have $\lambda_{\min}(G^{-1}\hat{\mathcal{G}}G^{-1}) \geq \lambda_{\min}(G^{-1}\mathcal{G}G^{-1}) - \|G^{-1}(\hat{\mathcal{G}} - \mathcal{G})G^{-1}\|$, which is bounded away from zero with probability tending to 1. It follows that $(G^{-1}\hat{\mathcal{G}}G^{-1})^{-1} = O_{\mathrm{P}}(1)$. Furthermore, we have

$$\Delta = (G^{-1}\hat{\mathcal{G}}G^{-1})^{-1} - (G^{-1}\mathcal{G}G^{-1})^{-1}$$
$$= -(G^{-1}\mathcal{G}G^{-1})^{-1}G^{-1}(\hat{\mathcal{G}} - \mathcal{G})G^{-1}(G^{-1}\hat{\mathcal{G}}G^{-1})^{-1} \longrightarrow 0$$

in probability. Therefore, we have

$$G(\hat{\Sigma}_{\mathrm{R}} - \Sigma_{\mathrm{R}})G = (G^{-1}\mathcal{G}G^{-1})^{-1}G^{-1}\mathcal{I}_1 G^{-1}\Delta + \Delta G^{-1}\mathcal{I}_1 G^{-1}(G^{-1}\hat{\mathcal{G}}G^{-1})^{-1}$$
$$+ (G^{-1}\hat{\mathcal{G}}G^{-1})^{-1}G^{-1}(\hat{\mathcal{I}}_1 - \mathcal{I}_1)G^{-1}(G^{-1}\hat{\mathcal{G}}G^{-1})^{-1} \longrightarrow 0$$

in probability, using the results previously proved.

8.4. *Proof of Theorem* 4.  We have

$$\hat{\chi}^2 = (K'\hat{\theta} - \varphi)'(K'\Sigma_{\mathrm{R}}K)^{-1}(K'\hat{\theta} - \varphi)$$
$$+ (K'\hat{\theta} - \varphi)'\{(K'\hat{\Sigma}_{\mathrm{R}}K)^{-1} - (K'\Sigma_{\mathrm{R}}K)^{-1}\}(K'\hat{\theta} - \varphi)$$
$$= \tilde{\chi}^2 + \Delta.$$

By (26), $\tilde{\chi}^2 \to \chi_r^2$, so it remains to show that $\Delta \to 0$ in probability.

According to the definition in Section 2.2 (first paragraph) and the conclusion of Theorem 2, for any $\eta > 0$ we have, with probability tending to 1 (hereafter w.p. $\to$ 1), $(1 - \eta)\Sigma_{\mathrm{R}} \leq \hat{\Sigma}_{\mathrm{R}} \leq (1 + \eta)\Sigma_{\mathrm{R}}$. It follows that w.p. $\to$ 1 $(1 - \eta)K'\Sigma_{\mathrm{R}}K \leq K'\hat{\Sigma}_{\mathrm{R}}K \leq (1 + \eta)K'\Sigma_{\mathrm{R}}K$ and hence $(1 + \eta)^{-1}(K'\Sigma_{\mathrm{R}}K)^{-1} \leq (K'\hat{\Sigma}_{\mathrm{R}}K)^{-1} \leq (1 - \eta)^{-1}(K'\Sigma_{\mathrm{R}}K)^{-1}$ (e.g., [20], Theorem A.52). Therefore, we have w.p. $\to$ 1,

$$(36) \quad (1 + \eta)^{-1}I_r \leq (K'\Sigma_{\mathrm{R}}K)^{1/2}(K'\hat{\Sigma}_{\mathrm{R}}K)^{-1}(K'\Sigma_{\mathrm{R}}K)^{1/2} \leq (1 - \eta)^{-1}I_r,$$

which implies $\|W - I_r\| \leq \{(1 - \eta)^{-1} - 1\} \vee \{1 - (1 + \eta)^{-1}\}$, where $W$ is the middle term in (36) and $a \vee b = \max(a, b)$. Since $\eta$ is arbitrary, this proves that $W \to I_r$ in probability. That $\Delta \to 0$ in probability then follows by observing

$$\Delta = (K'\hat{\theta} - \varphi)'(K'\Sigma_{\mathrm{R}}K)^{-1/2}(W - I_r)(K'\Sigma_{\mathrm{R}}K)^{-1/2}(K'\hat{\theta} - \varphi),$$

hence $\|\Delta\| \leq \|W - I_r\|(K'\hat{\theta} - \varphi)'(K'\Sigma_{\mathrm{R}}K)^{-1}(K'\hat{\theta} - \varphi) = \|W - I_r\|\tilde{\chi}^2$.

**Acknowledgments.**  The author is grateful to an Associate Editor and two referees for their constructive comments. Furthermore, the author wishes to thank Zhonghua Gu and Wen-Ying Feng for their computational support during the revision of the manuscript.

## REFERENCES

[1] ARVESEN, J. N. (1969). Jackknifing $U$-statistics. *Ann. Math. Statist.* **40** 2076–2100. MR0264805

[2] ARVESEN, J. N. and SCHMITZ, T. H. (1970). Robust procedures for variance component problems using the jackknife. *Biometrics* **26** 677–686.

[3] BATTESE, G. E., HARTER, R. M. and FULLER, W. A. (1988). An error-components model for prediction of county crop areas using survey and satellite data. *J. Amer. Statist. Assoc.* **83** 28–36.

[4] BIRKES, D. and WULFF, S. S. (2003). Existence of maximum likelihood estimates in normal variance-components models. *J. Statist. Plann. Inference* **113** 35–47. MR1963033

[5] DAS, K., JIANG, J. and RAO, J. N. K. (2004). Mean squared error of empirical predictor. *Ann. Statist.* **32** 818–840. MR2060179

[6] DEMIDENKO, E. and MASSAM, H. (1999). On the existence of the maximum likelihood estimate in variance components models. *Sankhyā Ser. A* **61** 431–443. MR1743550

[7] EFRON, B. and HINKLEY, D. V. (1978). Assessing the accuracy of the maximum likelihood estimator: Observed versus expected Fisher information (with discussion). *Biometrika* **65** 457–487. MR0521817

[8] GHOSH, M. and RAO, J. N. K. (1994). Small area estimation: An appraisal (with discussion). *Statist. Sci.* **9** 55–93. MR1278679

[9] HARTLEY, H. O. and RAO, J. N. K. (1967). Maximum likelihood estimation for the mixed analysis of variance model. *Biometrika* **54** 93–108. MR0216684

[10] HEYDE, C. C. (1994). A quasi-likelihood approach to the REML estimating equations. *Statist. Probab. Lett.* **21** 381–384. MR1325214

[11] HEYDE, C. C. (1997). *Quasi-Likelihood and Its Application*. Springer, New York. MR1461808

[12] JIANG, J. (1996). REML estimation: Asymptotic behavior and related topics. *Ann. Statist.* **24** 255–286. MR1389890

[13] JIANG, J. (1997). Wald consistency and the method of sieves in REML estimation. *Ann. Statist.* **25** 1781–1803. MR1463575

[14] JIANG, J. (1998). Asymptotic properties of the empirical BLUP and BLUE in mixed linear models. *Statist. Sinica* **8** 861–885. MR1651513

[15] JIANG, J. (2003). Empirical method of moments and its applications. *J. Statist. Plann. Inference* **115** 69–84. MR1972940

[16] JIANG, J. (2004). Dispersion matrix in balanced mixed ANOVA models. *Linear Algebra Appl.* **382** 211–219. MR2050107

[17] KHURI, A. I., MATHEW, T. and SINHA, B. K. (1998). *Statistical Tests for Mixed Linear Models*. Wiley, New York. MR1601351

[18] LAIRD, N. M. and WARE, J. H. (1982). Random effects models for longitudinal data. *Biometrics* **38** 963–974.

[19] LANGE, N. and RYAN, L. (1989). Assessing normality in random effects models. *Ann. Statist.* **17** 624–642. MR0994255

[20] RAO, C. R. and TOUTENBURG, H. (1995). *Linear Models*. Springer, Berlin. MR1354840

[21] RAYNER, J. C. W. and BEST, D. J. (1989). *Smooth Tests of Goodness of Fit*. Oxford Univ. Press, London. MR1029526

[22] RICHARDSON, A. M. and WELSH, A. H. (1994). Asymptotic properties of restricted maximum likelihood (REML) estimates for hierarchical mixed linear models. *Austral. J. Statist.* **36** 31–43. MR1309503

[23]  Searle, S. R., Casella, G. and McCulloch, C. E. (1992). *Variance Components.*
        Wiley, New York. MR1190470

Department of Statistics
University of California
One Shields Avenue
Davis, California 95616
USA
e-mail: jiang@wald.ucdavis.edu